\documentclass[12pt]{article}
\usepackage{amsthm,amsfonts,amsmath,amssymb}
\renewcommand{\proof}{\noindent{\it Sketch of proof~\,\,}}

\newcommand{\eq}{\begin{equation}}
\newcommand{\en}{\end{equation}}

\newcommand{\dd}{{\rm d }}
\newcommand{\giv}{\,|\,}

\newcommand{\darr}{\downarrow}

\newcommand{\Var}{{\rm Var}}

\newcommand{\prob}{\mathbb P}
\newcommand{\ex}{\mathbb E}

\newcommand{\Nat}{\Bbb N}

\newcommand{\ed}{ \stackrel{d}{=}}

\newcommand{\nut}{\tilde{\nu}}

\newcommand{\KK}{{\cal K}}

\newcommand{\up}{\uparrow}

\def\endpf{\hfill $\Box$ \vskip0.5cm}

\newtheorem{theorem}{\large Theorem}[section]

\newtheorem{definition}[theorem]{\large Definition}
\newtheorem{corollary}[theorem]{\large Corollary}
\newtheorem{lemma}[theorem]{\large Lemma}

\newtheorem{example}[theorem]{\large Example}

\begin{document}
\title{Regeneration in Random Combinatorial Structures}
\author{Alexander V. Gnedin\thanks{Postal address:
 Department of Mathematics, Utrecht University,
 Postbus 80010, 3508 TA Utrecht, The Netherlands. email: A.V.Gnedin@uu.nl}}
\date{}
\maketitle

\begin{abstract}
\noindent
Theory of Kingman's partition structures has two culminating points
\begin{itemize}
\item the general paintbox representation, relating finite partitions
 to hypothetical infinite populations 
via a natural sampling procedure,
\item a central example of the theory: the Ewens-Pitman two-parameter partitions. 
\end{itemize}
\noindent
In these notes we further develop the theory by 
\begin{itemize}
\item passing to structures enriched by the order on the collection of categories,
\item extending the class of tractable models by exploring the idea of regeneration,
\item analysing regenerative properties of the Ewens-Pitman partitions,
\item studying asymptotic features of  the regenerative compositions. 
\end{itemize}
\end{abstract}

\section{Preface}

The kind of discrete regenerative phenomenon discussed here
is present in the cycle patterns of random permutations.
 To describe this instance, first
recall that
every  permutation of $[n]:=\{1,\ldots,n\}$ is decomposable  in a product of disjoint cycles.
The cycle sizes  make up  a partition of $n$ into some number of positive integer parts. 
For instance, permutation 
$(1~3)(2)$ of the set $[3]$ corresponds to the partition of integer $3$ 
with  parts $2$ and $1$.
 Permutations of different degrees $n$ are connected in a natural way.
Starting with a permutation of $[n]$, a permutation of the smaller set $[n-1]$ is created 
by removing element $n$ from its cycle.  
This reduction is a surjective $n$-to-1 mapping. For instance, three permutations
$(1~3)(2),~(1)(2~3),~(1)(2)(3)$ are mapped  to $(1)(2)$.

\par Now suppose  the permutation is chosen uniformly at random from the set of all $n!$ permutations of  $[n]$. 
The collection of cycle-sizes is then a certain random partition $\pi_n$ of integer $n$.
By the $n$-to-$1$ property of the projection, 
the permutation reduced by element $n$ 
is the uniformly distributed permutation of $[n-1]$,
with the cycle partition $\pi_{n-1}$.
The transition from $\pi_n$ to $\pi_{n-1}$ is easy to describe directly, without reference to underlying permutations: 
choose a random part of $\pi_n$ by a size-biased pick, i.e.
with probability
proportional to the size of the part, and 
then reduce  the chosen part by $1$.
This transition rule   
suggests to view  
the random partitions with varying $n$ altogether as components
of an infinite 
{\it partition structure} $(\pi_n,~ n=1,2,\ldots)$.

\par  Apart from the consistency property inherent to any  partition structure, 
there is another recursive self-reproduction property of the partitions derived from the 
cycle patterns of 
uniform permutations.
Fix $n$ and suppose a part is chosen by a size-biased pick from  $\pi_n$ and  completely deleted.
Given the part was $m$, the partition reduced by this part will be a distributional copy of 
 $\pi_{n-m}$. 
In this sense
the partition structure $(\pi_n,~ n=1,2,\ldots)$  {\it regenerates}. 

\par  For large $n$, the size-biased pick  will  choose a part with about $nU$ elements, where $U$ is a 
random variable with 
uniform distribution on the unit interval. 
In the same way, the iterated deletion of parts by size-biased picking 
becomes similar
to the splitting of $[0,1]$ at points representable via products of independent uniform variables.
 The latter is a special case  
of the multiplicative renewal process often called {\it stick-breaking}. 

\par In these notes we consider sequences of partitions and ordered partitions which are consistent in the same 
sense as the cycle patterns of permutations for various $n$. 
In contrast to that, the assumption about the regeneration property of such structures will be fairly general. 
The connection between  combinatorial partitions and  splittings of the unit interval 
is central in the theory and will be analysed in detail
in  the  general context of regenerative structures.

\section{The paintbox and the two-parameter family}

A {\it composition} of integer $n$ is an ordered sequence $\lambda^\circ=(\lambda_1,\ldots,\lambda_k)$ of positive integer parts with 
sum $|\lambda^\circ|:=\sum_j\lambda_j=n$.
We shall think of composition as a model of occupancy, meaning $n$ `balls' separated by `walls' into some number of nonempty `boxes', 
like in this diagram 
$${\bf |}\bullet\bullet\bullet{\bf |}\bullet{\bf |}\bullet\bullet\,\,{\bf |}$$
representing composition $(3,1,2)$. A wall ${\bf |}$ is either placed between two consequitive $\bullet$'s or not, hence there are $2^{n-1}$ compositions of 
$n$. 
Sometimes we shall also use encoding the compositions into binary sequences, in which 
a $1$ followed by some $m-1$  zeroes corresponds to part $m$, like the code
$100110$ for composition $(3,1,2)$,

\par 
A related labeled object
is an {\it ordered partition} of the set $[n]:=\{1,\ldots,n\}$, which may be obtained by some enumeration of the 
balls by integers $1,\ldots,n$, like
$$ |\stackrel{2}{\bullet}\,\,\stackrel{4}{\bullet}\,\,\stackrel{5}{\bullet}|\stackrel{3}{\bullet}|\stackrel{1}{\bullet}\,\,\stackrel{6}{\bullet}|$$
(the ordering of balls within a box is not important). 
The number of such labelings, that is the number of ordered set partitions with {\it shape} $(\lambda_1,\ldots,\lambda_k)$, is equal to 
the multinomial coefficient 
$$f^\circ(\lambda_1,\ldots,\lambda_k):={n!\over \lambda_1!\cdots\lambda_k!}.$$
Throughout, symbol $^\circ$ will denote 
a function of composition, also when the function is not sensitive to the permutation of parts.

\par Discarding the order of parts in a composition $(\lambda_1,\ldots,\lambda_k)$ yields a {\it partition} of integer 
$|\lambda|$, usually written as a ranked sequence of nondecreasing parts.
For instance, the ranking maps compositions $(3,1,2)$ and $(1,3,2)$ to the same partition
$(3,2,1)^\downarrow$,  where $\downarrow$ will be both used to denote the operation of ranking and to
indicate that the arrangement of parts in sequence is immaterial. 
Sometimes we use notation like $2\in (4,2,2,1)^\darr$ to say that $2$ is a part of partition. 
The number of partitions of the set $[n]$ with the same shape 
$\lambda^\downarrow=(\lambda_1,\ldots,\lambda_k)^\downarrow$ is equal to
$$f(\lambda^\downarrow):=n!\prod_{r=1}^n {1\over (r!)^{k_r}k_r!}~$$
where $k_r=\#\{j:\lambda_j=r\}$ is the number of parts of $\lambda^\downarrow$ of size $r$.

\par A random composition/partition of $n$ is simply a random variable with values in the  finite set 
of compositions/partitions of $n$. One statistical context where these combinatorial objects appear is the 
{\it species sampling problem}. Imagine an alien who has no idea of the mammals.  
Suppose the first six mammals she observes are tiger, giraffe, elephant, elephant, elephant and  giraffe, appearing in this 
sequence.
Most frequent -- three of these -- have long trunks, two are distinctively taller than the others, 
and one is striped.
She records this as partition $(3,2,1)^\downarrow$ into three distinct species. 
Composition $(1,3,2)$ could appear 
as the record of species abundance by more delicate classification  
according to typical height, from the lowest to the tallest\footnote{If  her guidebook would describe 
four species, e.g. these three and the cows,  her records would be
$(3,2,1,0)^\downarrow,(1,0,3,2)$ 
(weak partitions, respectively, weak compositions), 
but we assumed that she knew apriori really 
nothing of the mammals.}. 
Enumerating the animals in the order of observation gives a
labeled object, a
 partition/ordered-partition of the set 
$[6]=\{1,\ldots,6\}$.

\par There are many ways to introduce random partitions or compositions.
The method adopted here is
intrinsically related to the species sampling problem.
This is the following ordered version of Kingman's paintbox
(see  \cite{Bertoin},  \cite{GnedinRCS},  \cite{PitmanBernoulli}).

\vskip0.2cm
\noindent
{\bf Ordered paintbox}
{\rm 
Let $\cal R$ be a random closed subset of $[0,1]$.
The complement open set ${\cal R}^c:=(0,1)\setminus {\cal R}$ has a canonical representation
as a disjoint union of countably many open  interval components, which we shall call the {\it gaps} of $\cal R$.  
Independently of $\cal R$, sample points 
$U_1, U_2,\ldots$ from the uniform distribution on $[0,1]$ and group the points
in clusters by the rule:
 $U_i, U_j$ belong to the same cluster if they
hit the same gap of $\cal R$.
If  $U_i$ falls in 
${\cal R}$ let $U_i$ be a singleton.
For each $n$, count the representatives of clusters among $U_1,\ldots,U_n$ and define $\varkappa_n$,
a random  composition  of integer $n$, to be the 
record of positive counts in the  left-to-right order of the gaps.}

\vskip0.2cm
\noindent
For instance,  $\varkappa_n$ assumes the value $(3,1,2)$ if, in the left-to-right order,
there is a gap hit by three 
points out of $U_1,\ldots,U_6$,
a singleton cluster resulting from either some 
gap or from some $U_j\in{\cal R}$, and a gap hit by two of $U_1,\ldots,U_6$.

\vskip0.2cm
\noindent
In the {\it proper} case  ${\cal R}$ has Lebesgue measure zero almost surely, hence
 $U_j\in {\cal R}$ occurs only with probability zero.
We may  think then of points of $\cal R$ as possible locations of walls $|$ and of the points of $[0,1]$ 
as possible locations of balls $\bullet$. 
In a particular realisation, the balls appear at locations
$U_j$, and the walls bound the gaps hit by at least one ball.
In the {\it improper} case, $\cal R$ may have positive measure with nonzero probability. 
If $U_j\in {\cal R}$
we can imagine a box with walls coming so close together that no further ball will fit in this box,
so $U_j$ will forever remain a singleton, no matter how many balls are added.

\par Sometimes we shall identify $\cal R$ with the splitting of $[0,1]$ it induces, and just call $\cal R$ itself the paintbox.
Molchanov \cite{Molchanov} gives  extensive exposition of the theory of random sets,
although an intuitive idea will 
suffice from most of our purposes.
This can be a set of some fixed cardinality, e.g.  
splittings of $[0,1]$ following a Dirichlet distribution (see \cite{PD}, \cite{Huilletmart}), 
or complicated random Cantor-type sets like the set of zeroes of the Brownian motion.
It will be also convenient to make no difference between   
two closed subsets of $[0,1]$ if they
only differ by endpoints $0$ or $1$. 
If $1$ or $0$ is not  accumulation point for $\cal R$, the gap adjacent to the boundary 
 will be called  right or left {\it meander}.

\par The paintbox with random $\cal R$ is a kind of canonical representation 
of `nonparametric priors' in the species sampling problem. 
View  ${\mathbb R}$ as an ordered space of 
distinct {\it types}. 
Originally, by Kingman \cite{KingmanRPS}, 
the types  were colours making up a paintbox.
Consider a random probability measure
$F$ on reals as a model of infinite ordered population. 
Let $\xi_1,\xi_2,\ldots$ be a sample from $F$, which means that 
conditionally given $F$, the $\xi_j$'s are i.i.d. with distribution $F$.
An ordered  partition of the sample
is defined by grouping $j$'s with the same value of $\xi_j$,
with the order on the groups
maintained by increase of the values.
The case of diffuse (nonatomic) $F$ is trivial --
then ties among $\xi_j$'s have probability zero and the partition 
has only singletons, so the substantial case is  
$F$ with atoms, when  the partition will have nontrivial blocks.
The same ordered partition is induced by any other distribution obtained from $F$ by a suitable 
monotonic transformation, which may be random.  
To achieve the uniqueness, view $F$ as a random distribution function 
and observe that $\xi_i\leq \xi_j$ iff $F(\xi_i)\leq F(\xi_j)$.
Conditioning on $F$ and applying the quantile transform $y\to F(y)$ to the sample produces  
another sample $\tilde{\xi}_1,\tilde{\xi}_2,\ldots$  from the transformed distribution
 $\tilde{F}$ supported by $[0,1]$.
In the diffuse case, $F$ is well known to be  the uniform distribution, and in general
the distribution function $\tilde{F}$ is of special kind: it
satisfies $F(x)\leq x$ for $x\in [0,1]$ and $\tilde{F}(x)=x$ $\tilde{F}$-a.s.
Moreover, 
each jump location of $\tilde{F}$ is preceded by a flat (where $\tilde{F}$ is constant), whose length is equal to the size of the 
jump. The latter implies that the composition derived from $\tilde{F}$ by grouping equal $\tilde{\xi}_j$'s 
 in clusters is the same as the composition obtained via the paintbox construction from ${\cal R}={\rm support}(\tilde{F})$.
The identification with the paintbox construction  can be  shown more directly, i.e.
without appealing to $\tilde{F}$,
by taking for $\cal R$ the range of the random function $F$
(note that ${\rm support}(\tilde{F})$ with $0$ attached to it coincides with the range of $F$).

\vskip0.2cm
Note further important features inherent to the  paintbox construction:
\begin{itemize}
\item 
The unlabeled object, $\varkappa_n$, is determined by $\cal R$ and the uniform order statistics
$U_{n:1}<\ldots<U_{n:n}$, i.e. the ranks of $U_1,\ldots,U_n$ appear as random labels and do not matter.
\item Attaching label $j$ to the ball corresponding to $U_j$, we obtain, for each $n$, an  ordered partition ${\rm K}_n$
of the set $[n]$, with shape $\varkappa_n$. 
This ordered partition is {\it exchangeable}, meaning that a permutation of the labels does not change
the distribution of ${\rm K}_n$, thus
all ordered partitions of $[n]$ with the same shape have the same 
probability.

\item The ordered partitions ${\rm K}_n$ are {\it consistent} as $n$ varies.
Removing ball $n$ (and deleting an empty box in case one is created) reduces ${\rm K}_n$ to ${\rm K}_{n-1}$.
The infinite sequence ${\rm K}=({\rm K}_n)$ of consistent ordered partitions of
$[1], [2],\ldots$ 
defines therefore an exchangeable ordered partition of the infinite set $\Nat$ into some collection of nonempty blocks. 

\end{itemize}

\noindent
Translating the consistency in terms of  compositions $\varkappa_n$ we arrive at

\begin{definition}\label{def1} {\rm A sequence $\varkappa=(\varkappa_n)$ of random compositions of $n=1,2,\ldots$ is called a 
{\it composition structure}
if these are {\it sampling consistent}: for each $n>1$, conditionally given $\varkappa_n=(\lambda_1,\ldots,\lambda_k)$
the composition $\varkappa_{n-1}$ has the same distribution as the composition obtained by reducing by $1$ each part $\lambda_j$ with
probability $\lambda_j/n$.
}
\end{definition}

\par A {\it size-biased~} part of  composition $\lambda^\circ$ is a random part which coincides with every part $\lambda_j$ 
with probability $\lambda_j/|\lambda^\circ|$. 
A size-biased part of a random composition $\varkappa_n$ is defined conditionally on the value $\varkappa_n=\lambda^\circ$.
The sampling consistency condition amounts to the transition 
from $\varkappa_n$ to $\varkappa_{n-1}$  by 
reducing  a size-biased part. 
This special reduction rule in Definition \ref{def1} 
is a trace of the exchangeability in ${\rm K}_n$ that remains when the labels are erased: 
indeed, given the sizes of the blocks,
the ball with label $n$  belongs to a particular block of size $\lambda_j$ with probability $\lambda_j/n$.  
\par Keep in mind that
 the consistency of ordered set partitions ${\rm K}_n$ is understood in the {\it strong~} sense, as a property of  
random objects defined on the same probability space,
while Definition \ref{def1} only requires  {\it weak} consistency in terms of the distributions of $\varkappa_n$'s.
By the measure extension theorem, however, the correspondence between (the laws of) exchangeable ordered partitions of $\Nat$
and composition structures is one-to-one, and any composition structure can be realised through an exchangeable ordered partition of 
$\Nat$. 
In view of this correspondence, dealing with labeled or unlabeled objects is just the matter of convenience, and we
shall freely switch from one model to another.

\par A central  result about the general composition structures says that these can be uniquely 
represented by a paintbox \cite{GnedinRCS}.
This extends Kingman's \cite{KingmanRPS} representation of partition structures.

\begin{theorem}\label{Thm1} 
For every composition structure $\varkappa=(\varkappa_n)$ there exists a unique distribution for a random closed set $\cal R$ which 
by means of the paintbox construction
yields, for each $n$,
a distributional copy of  $\varkappa_n$. 
\end{theorem}
\proof
 The line of the proof is analogous to modern proofs of 
de Finetti's theorem which asserts that a sequence of exchangeable random variables is conditionally  i.i.d.
given the limiting empirical distribution of the sequence (see Aldous \cite{Aldous}). 
To this end, we need to make  the concept of a random closed set precise.
One way to do this is to topologise 
the space of closed subsets of $[0,1]$ 
by means of
the Hausdorff distance. Recall that for
$R_1,R_2\subset[0,1]$ (with boundary points $0,1$ adjoined to the sets) the distance is 
equal to the smallest $\epsilon$ such that the $\epsilon$-inflation of $R_1$
covers $R_2$ and the same holds with the roles swapped, so the distance is small when
the sizes and positions of a few biggest gaps are approximately the same for both sets.
Realise all $\varkappa_n$'s on the same probability space through some exchangeable 
${\rm K}$. Encode each composition $(\lambda_1,\ldots,\lambda_k)$ into 
a finite set $\{0,\Lambda_1/n,\ldots,\Lambda_{k-1}/n,1\}$ where $\Lambda_j=\lambda_1+\ldots+\lambda_j$.
This maps $\varkappa_n$ to a finite random set ${\cal R}_n\subset [0,1]$. 
By a martingale argument it is shown that  
the law of the large numbers holds: 
as $n\to\infty$ 
the sets 
${\cal R}_n$ converge almost surely to a random closed set $\cal R$. 
The limit $\cal R$ is shown to direct the paintbox representation of $\varkappa$.
\endpf
\vskip0.3cm


There are various  equivalent formulations of the result in terms of (i) the exchangeable quasi-orders on $\Nat$
(in the spirit of \cite{Jacka}),
(ii) the entrance Martin boundary for the time-reversed Markov chain $(\varkappa_n, n=\ldots,2,1)$, (iii) certain
functionals on the infinite-dimensional algebra of quasisymmetric functions
\cite{GO}.
\vskip.2cm

\par 
We define the {\it composition probability function} (CPF for shorthand) 
$p^\circ$  of a composition structure $\varkappa$ 
as
$$p^\circ(\lambda^\circ):={\mathbb P}(\varkappa_n=\lambda^\circ),~~~|\lambda^\circ|=n,~n=1,2,\ldots $$
For fixed $|\lambda^\circ|=n$ this is the distribution of $\varkappa_n$.
To avoid 
confusion with the distribution of ${\rm K}_n$ we stress that
the probability of any {\it particular} value of the set partition ${\rm K}_n$ with
shape $\lambda^\circ$ is equal to
 $p^\circ(\lambda^\circ)/f^\circ(\lambda^\circ)$. 
Sampling consistency translates as a backward recursion
\begin{equation}\label{back}
p^\circ(\lambda^\circ)=\sum_{\mu^\circ} c(\lambda^\circ,\mu^\circ) p^\circ(\mu^\circ),
\end{equation}
where $\mu^\circ$  runs over all shapes of  extensions of any fixed  ordered partition of $[n]$ with shape $\lambda^\circ$ 
to some ordered partition of $[n+1]$. For instance, taking $\lambda^\circ=(2,3)$, 
$\mu^\circ$ assumes the values $(1,2,3),(2,1,3),(2,3,1),(3,3),(2,4)$.
The coefficient
$c(\lambda^\circ,\mu^\circ)$ is the probability to obtain $\lambda^\circ$ from $\mu^\circ$ by reducing 
a size-biased part of $\mu^\circ$.

\par For fixed $n$, if $p^\circ(\lambda^\circ)$ is known for compositions $\lambda^\circ$ with $|\lambda^\circ|=n$,
then solving (\ref{back}) backwards gives the values of CPF for all compositions with $|\lambda^\circ|\leq n$.
By linearity of the recursion, every such partial solution, with $n'\leq n$, is a convex combination of 
$2^{n-1}$ solutions obtained by taking delta measures on the level $n$.
Similarly, without restricting $n$,
the set of CPF's is 
convex and compact in the weak topology
of functions on a countable set; 
 this convex set has the property
of uniqueneess of barycentric decomposition in terms of extreme elements
(Choquet simplex). 
The extreme CPF's are precisely those 
derived from nonrandom paintboxes.
The correspondence between 
extreme solutions and  closed subsets of $[0,1]$ is a homeomorphism, which extends to 
the homemorphism between all CPF's and distributions for random closed $\cal R$.

\vskip0.2cm
\par Discarding the order of parts in each $\varkappa_n$ we  obtain Kingman's {\it partition structure}
$\pi=(\pi_n)$ with $\pi_n=\varkappa_n^\downarrow$. 
Partition structures satisfy the same sampling consistency condition as in Definition \ref{def1}.
The corresponding labeled object
is an {\it exchangeable partition} $\Pi=(\Pi_n)$ of the infinite set $\Nat$.
The law of large numbers for partition structures says that,
as $n\to\infty$, the vector $n^{-1}\pi_n$ 
padded by infinitely many zeroes
converges (weakly for $\pi_n$, strongly for $\Pi_n$)
to a random element 
 $\cal S$ of the 
infinite-dimensional simplex
$$\nabla=\{(s_i):\,\,s_1\geq s_2\ldots  \geq 0, \,\sum_is_i\leq 1\},$$
so the 
components of $\cal S$ are the asymptotic {\it frequencies} of the ranked parts of $\varkappa_n$. 
 The {\it partition probability function} (PPF) 
$$p(\lambda^\darr):=\prob(\pi_n=\lambda^\darr),~~~~|\lambda^\darr|=n,~n=1,2,\ldots,$$
specifies  distributions of $\pi_n$'s and satisfies a recurrence analogous to (\ref{back}).
The correspondence between PPF's and distributions for {\it unordered paintbox} $\cal S$
is bijective. Note that the possibility of strict inequality  $\sum_j s_j<1$ occurs in the improper case,
where the diffuse mass $1-\sum_j s_j$, sometimes  also called dust \cite{Bertoin}, 
is equal to the cumulative frequency of singleton blocks of $\Pi$ given ${\cal S}=(s_j)$.

\vskip0.2cm
\par Discarding order is a relatively easy operation.
In terms of ordered and unordered paintboxes  $\cal R$ and $\cal S$  the connection is expressed by the formula
\eq\label{SU}
{\cal S}=({\cal R}^c)^\downarrow,
\en
where the ranking $\darr$ means 
that the gap-sizes of $\cal R$ are recorded in nonincreasing order.
The operation $\downarrow$ 
is a continuous mapping from the space
of closed subsets of $[0,1]$ to $\nabla$. 
In terms of distributions, passing from CPF to PPF is expresses by
 the symmetrisation formula
\eq\label{symm1}
p(\lambda^\downarrow)=\sum_\sigma p^\circ(\lambda^\sigma),
\en
where $\lambda^\sigma$ runs over all distinct arrangements of parts of $\lambda^\downarrow$ in a composition
(e.g. for partition $(2,1,1)^\downarrow$ there are three such compositions $(2,1,1),(1,2,1),(1,1,2)$).

\par  In the other direction, there is one universal way to introduce the order.
With every partition structure one can accosiate a unique {\it symmetric} composition structure, for which
any of the following three equivalent conditions holds:

\begin{itemize}
\item[(i)] 
all terms in the RHS of (\ref{symm1}) are equal,
\item[(ii)] conditionally given $\Pi_n$ with $k$ blocks, any arrangement of the blocks in ${\rm K}_n$ has the same probability,
\item[(iii)] the gaps of $\cal R$ appear in the exchangeable random order.
\end{itemize}
The last property (iii)  means that,  
 conditionally given ${\cal S}=(s_j)$ with $s_k>0$, 
every relative order of the first $k$ largest gaps (labeled by $[k]$) 
of sizes $s_1,\ldots,s_k$ has probability $1/k!$.  This rule defines $\cal R$ unambiguously 
in the proper case, and extension to the
improper case follows  by continuity. A simple example of symmetric $\cal R$ is associated with splitting $[0,1]$ 
according to the symmetric Dirichlet distribution on a finite-dimensional simplex.

\par Beside from the symmetric composition structure, there are many other composition structures 
associated with a given partition structure.
Understanding the connection in the direction from unordered to ordered structures
 is a difficult {\it problem of arrangement}.
To outline some facets of the problem,
suppose we have a rule to compute $p^\circ$ from $p$, 
how can we pass then from $\cal S$ to $\cal R$?
Specifically,
given ${\cal S}=(s_j)$,
in which order  the intervals of sizes $s_1,s_2,\ldots$ should be arranged in an open set?
Other way round, 
suppose we have a formula for $p$ and know that (\ref{SU}) is true,
 how then 
can we compute the probability  that given
$\pi_5=(3,2)^\downarrow$ 
the parts appear in the composition as $(2,3)$?
Most questions like that 
cannot have universal  answers, because 
random sets and random series are objects of high complexity,
and the paintbox correspondence cannot be expressed
by simple formulas. 

\vskip0.2cm
\noindent
{\bf Ewens-Pitman partition structures} In the theory of partition structures 
and partition-valued processes of fragmentation and coagulation \cite{Bertoin}
a major role is played by 
the Ewens-Pitman two-parameter family of partitions, with PPF
\eq\label{2par}
p_{\alpha,\theta}
(\lambda^\downarrow)=f(\lambda^\downarrow)
{\prod_{i=1}^{k-1} (\theta+\alpha i)\over (1+\theta)_{n-1}}\prod_{j=1}^k (1-\alpha)_{\lambda_j-1},
~~~~~~\lambda^\downarrow=
(\lambda_1,\ldots,\lambda_k)^\downarrow,
\en
where and henceforth $(z)_n:=z (z+1)\cdots (z+n-1)$ 
is a rising factorial. 
The principal range of the parameters is 
\eq\label{2par-range}
\{(\alpha, \theta):  0\leq\alpha< 1, \theta>-\alpha\}\cup \{(\alpha, \theta):\alpha<0, -\theta/\alpha\in \Nat\},
\en
and there are also a few  
degenerate  boundary cases defined by continuity.
\par One of many remarkable features of these partitions 
is the sequential device
for generating 
the corresponding exchangeable partition $\Pi=(\Pi_n)$. Start with the one-element 
partition $\Pi_1$. Inductively,
Suppose $\Pi_n$ has been constructed then,
given that the shape of $\Pi_n$ is $(\lambda_1,\ldots,\lambda_k)^\downarrow$, the ball $n+1$ is placed in the existing box $i$ with probability
$(\lambda_i-\alpha)/(n+\theta)$ for $i=1,\ldots,k$, and  starts a new box with probability $(\theta+k\alpha)/(n+\theta)$.
In the Dubins-Pitman  interpretation as a
`Chinese restaurant process', the  
 balls correspond to customers arriving in 
the restaurant, 
and boxes  are  circular tables.
With account of the circular ordering of customers at each occupied table, 
and subject to uniform random placement at each particular table,
the process also defines a consistent sequence of 
random permutations for $n=1,2,\ldots$; with uniform distributions in the case $(\alpha,\theta)=(0,1)$.

\par  The two-parameter family has numerous connections to basic types of random processes like 
the Poisson process and the Brownian motion, see Pitman's lecture notes \cite{CSP} for a summary.   
It also provides an exciting framework for the problem of arrangement.

\section{Regenerative composition structures}

Every $\varkappa_n$ 
in a composition structure  may be regarded as  a reduced copy of $\varkappa_{n+1}$.
We  complement this now
by another type of self-reproduction property, related to the reduction by a whole box.

\begin{definition}
\label{regendef}
{\em A composition structure $\varkappa=({\varkappa}_{n})$ is called {\em regenerative}
if for all $n>m\geq 1$, the following deletion property holds.
If the first part of ${\varkappa}_{n}$ is deleted and
conditionally given  this part   is $m$,
the remaining composition of $n-m$ is distributed like ${\varkappa}_{n-m}$.
}
\end{definition}
\noindent

Denote $F_n$ 
the first part of $\varkappa_n$ and consider its distribution  
$$q(n:m):=\sum_{|\lambda^\circ|=n,\,\lambda_1=m} p^\circ(\lambda^\circ).$$
It follows immediately from the definition that 
$\varkappa$ 
is regenerative 
iff the CPF has the product form
\begin{equation}
\label{produ}
p^\circ(\lambda_1,\ldots,\lambda_k)= \prod_{j=1}^{k} q(\Lambda_j:\lambda_j),
\end{equation}
where $\Lambda_j=\lambda_j+\ldots+\lambda_k$ for $1\leq j\leq k$.

\par  For each $n$, the formula identifies  $\varkappa_n$ with the sequence 
of decrements of a decreasing Markov chain $Q_n^\darr=(Q_n^\darr(t),~t=0,1,\ldots)$ on $0,\ldots,n$.
The chain starts at $n,$ terminates at $0$, and jumps from $n'\leq n$ to $n'-m$ with probability $q(n':m)$.
The binary code of $\varkappa_n$ is obtained by writing 
$1$'s in positions $n-Q_n^\darr(t)+1,\,t=0,1,\ldots$, and writing $0$' is all other positions,
with the convention that the last $1$ in position $n+1$ is not included in the code.
In view of this interpretation, we 
call $q=(q(n:m), 1\leq m\leq n,~n\in \Nat)$ the {\it decrement matrix} of  $\varkappa$.
Since $p^\circ$ is computable from $q$, the decrement matrix determines completely the distributions of $\varkappa_n$'s
and the distribution of the associated exchangeable ordered partition $\rm K$.

\par For  a  given regenerative $\varkappa$ let $\pi=(\pi_n)$, with 
$\pi_n=\varkappa_n^\darr$, be the related partition structure.
Think of $\varkappa_n$ as an arrangement of parts of $\pi_n$ in some order.
For partition $\lambda^\darr$ of $n$ and each $m\in\lambda^\darr$ 
define the {\it deletion kernel}
$$d(\lambda^\darr,m)=\prob(F_n=m\giv \pi_n=\lambda^\darr),$$ 
which specifies the conditional probability, given the unordered multiset of parts,
 to place a part of size $m$ in the first position in  $\varkappa_n$
(so $d(\lambda^\darr,m)=0$ if $m\notin \lambda$).
The deletion property of $\varkappa$ implies that the PPF of $\pi$ satisfies the identity
\eq\label{pd}
p(\lambda^\darr)d(\lambda^\darr,m)=q(n:m)p(\lambda^\darr\setminus \{m\}),
\en
where $q(n:\cdot)$, the distribution of $F_n$, may be written in terms of the deletion kernel as
\eq\label{qq}
q(n:m)=\sum_{\{\lambda^\darr:\,|\lambda^\darr|=n,\, \,m\in \lambda^\darr\}}
d(\lambda^\darr,\,m)p(\lambda^\darr).
\en

\par Intuitively, the deletion kernel is a stochastic algorithm of choosing a part of partition $\pi_n$ to place it in the first position
of composition $\varkappa_n$. Iterated choices arrange all parts of each $\pi_n$ in $\varkappa_n$, hence the 
deletion kernel may be used to describe the arrangement on the 
level of finite partitions.
The partition structure $\pi$ inherits from $\varkappa$ the property of invariance 
under deletion of a part chosen by {\it some} random rule, 
expressed formally as  (\ref{pd}) and (\ref{qq}).
This is, of course, a subtle property when compared with more obvious invariance of $\varkappa$
under the first-part deletion, as specified in  Definition \ref{regendef}.


\subsection{Compositions derived from stick-breaking}

Exploiting the paintbox construction  
we shall give  a large family of examples of regenerative composition structures.
The method  is called
stick-breaking, and it is also known 
under many other names like e.g.
residual allocation model or, deeper in history, 
random alms \cite{Halmos}.

Let $(W_i)$ be independent copies of some random variable
$W$ with range $0 < W \le 1$. 
A value of $W$ is chosen, and the unit stick $[0,1]$ is broken at location $W$ in two pieces, then the left 
 piece of size $W$ is frozen, and the  right piece of size $1-W$
is broken again in proportions determined by another copy of $W$, and so on ad infinitum.
The locations of  breaks make up a random set $\cal R$ with points
\eq
\label{sbreak}
Y_k= 1-\prod_{i=1}^k (1-W_i), ~~k=1,2,\ldots, 
\en
so the gaps are ${\cal R}^c=\cup_{k=0}^\infty (Y_k,Y_{k+1})$.
The cardinality of $\cal R$ is  finite if ${\prob}(W=1) >0$,
but  otherwise infinite, with points accumulating only at the  right endpoint of the unit interval.
 By the i.i.d. property of the proportions,  the part of $\cal R$
to the right of $Y_1$ 
is a scaled copy of the whole set,
\eq\label{regR}
{({\cal R}\cap [Y_1,1])-Y_1\over 1-Y_1}\ed{\cal R},
\en
and this re-scaled part of $\cal R$ is independent of  $Y_1$.

\par Suppose a composition structure $\varkappa$ is derived from the paintbox ${\cal R}=\{Y_j,~j=0,1,\ldots\}$.
If $(0,Y_1)$ contains at least one of the first $n$ uniform points $U_j$, 
then  the first part of the composition $\varkappa_n$ is equal to 
the number of uniforms 
hitting this interval. Otherwise, conditionally given $Y_1$, 
the sample comes from the uniform  distribution on $[Y_1,1]$. 
Together with the property (\ref{regR})
of $\cal R$ this implies 
$$q(n:m)={n\choose m}{\ex \,}\left( W^m (1-W)^{n-m}\right) +{\ex \,}(1-W)^n\,q(n:m),$$
whence the law of the first part of $\varkappa_n$ is

\begin{equation}\label{stbr}
q(n:m)=\frac{{n\choose m}{\ex \,}\left(W^m (1-W)^{n-m}\right)}{{\ex \,}\left(1-(1-W)^n\right)} 
\qquad m=1,\ldots,n.
\end{equation}
which is a mixture of binomial distributions conditioned on a positive value. 
The key property (\ref{regR})  we exploited can be generalised for every $Y_k\in {\cal R}$, 
from which iterating the argument we obtain the product formula (\ref{produ}).

\par Concrete examples are obtained by choosing a distribution for $W$.
 For instance, taking delta measure $\delta_x$ with some $x\in (0,1)$ yields
 ${\cal R}=\{1-(1-x)^k, ~k=0,\ldots,\infty\}$,
which induces the same composition structure as the one
associated with 
sampling from  the geometric distribution on the set of integers. 
This composition structure was studied in many contexts, 
inluding theory of records and random search algorithms.

\par Expectations involved in (\ref{stbr}) may be computed explicitly only in some cases, e.g. 
for $W$ with polynomial density,  
but even then the product formula (\ref{produ}) rarely simplifies.

\vskip0.2cm
\noindent
{\bf Example} Here is an  example of a relatively simple decrement matrix. 
Taking $W$ with the general two-parameter beta density 
\eq\label{beta}
\nu(\dd x)={ x^{\gamma-1}(1-x)^{\theta-1}\dd x \over {\rm B}(\gamma,\theta)},
~~~(\gamma,\theta>0)
\en
 we arrive at
\eq\label{2parbeta}
q(n:m)={n\choose m} {(\gamma)_m(\theta)_{n-m}\over (\gamma+\theta)_n-(\theta)_n}.
\en
The product formula (\ref{produ}) simplifies  moderately  for general integer $\gamma$ \cite{TSF}, and massively
in the following case $\gamma=1$.

\vskip0.2cm
\noindent
 {\bf Regenerative composition structures associated with Ewens' partitions} 
Now suppose  $W$ has a beta$(1,\theta)$ density 
$$\nu(\dd x)=\theta(1-x)^{\theta-1}\dd x, ~~x\in (0,1).$$
Evaluating beta integrals in (\ref{stbr})  we find the decrement matrix
\begin{equation}\label{ESF-q}
q(n:m)={n\choose m}\frac{(\theta)_{n-m}\, m!}{(\theta+1)_{n-1}\,n}\,,
\end{equation}
and
massive cancellation in (\ref{produ}) gives the CPF
 \eq\label{ESF}
p^\circ_{0,\theta}(\lambda_1,\ldots,\lambda_k)={\theta^k n!\over (\theta)_n}\prod_{j=1}^k {1\over\Lambda_j},
\en
with $\Lambda_j=\lambda_j+\ldots+\lambda_k$.
Symmetrisation (\ref{symm1}) gives
the PPF known as the  {\it Ewens sampling formula} (ESF)
\begin{equation}\label{ESFpart}
p_{0,\theta}(\lambda^\downarrow)= f(\lambda^\downarrow)  {\theta^k\over (\theta)_n},
\end{equation}
which is a special case of (\ref{2par}).
Recall that 
the combinatorial factor is  the number of set partitions of $[n]$ with given shape.
The  range of parameter is  $\theta\in[0,\infty]$, with the boundary cases defined by continuity.

\par  For $\theta=1$, the distribution of $W$ is uniform$[0,1]$ and
$q(n:m)=n^{-1}$ is a discrete uniform distribution for each $n$;
the associated partition $\pi_n$ is the same as 
the cycle partition of a uniform random permutation of $[n]$.
For general $\theta$, the ESF corresponds to a biased permutation, which 
for each $n$ takes a particular value with probability $\theta^{\#{\rm cycles}}/(\theta)_n$.

\par We shall call $\varkappa$ with CPF  (\ref{ESF}) 
Ewens' regenerative composition structure.
The problem of arrangement has in this case
a  simple explicit solution.
For partition $(\lambda_1,\ldots,\lambda_k)^\darr$ the {\it size-biased permutation} is the random arrangement of parts
obtained by the iterated size-biased picking  without replacement.
For nonnegative $(s_j)\in\nabla$ with $\sum_j s_j=1$ 
define a  size-biased permutation in a similar way:
 a generic term $s_j$ is placed in position $1$ with probability proportional to $s_j$, then 
another term is chosen by a size-biased pick from the remaining terms and placed in position 2, etc.
The resulting random sequence
is then  in the {\it size-biased order}, hence the  distribution of the sequence is invariant under the size-biased permutation
\footnote{To define a size-biased arrangement in the improper case $\sum_j s_j\leq 1$
consider any closed set $R\subset[0,1]$ with gap-sizes $(s_j)$.
Sample uniformly balls $U_j$ and record  the gap-sizes by increase of the minimal
labels of balls, with understanding the points of $R$ as zero-size gaps.}.

\begin{theorem}\label{ESForder} 
Ewens' composition structure {\rm (\ref{ESF})} has  parts in the size-biased order, for every $n$.
Conversely, if a regenerative composition structure has parts in the size-biased order, then
its CPF is {\rm (\ref{ESF})} for some $\theta\in [0,\infty]$.  
\end{theorem}
\noindent
The paintbox also has a similar property: 
the intervals $(Y_{j},Y_{j+1})$ are in the size-biased order.
The law of 
frequencies $\cal S$ 
is known as Poisson-Dirichlet distribution.
The law of the gap-sizes $(Y_1-Y_0, Y_2-Y_1,\ldots)$ is called the GEM distribution.

\vskip0.2cm
\noindent
{\bf Remark} For set partitions, size-biased ordering 
is sometimes understood as 
the arranging of  blocks of $\Pi_n$ by increase of their minimal elements
(other often used names:  age ordering, sampling ordering).
This creates an ordered partition for each $n$, but this ordered partion is not exchangeable, 
since e.g. element $1$ is  always 
in the first block. In Ewens' case, but not in general,  the unlabeled compositions associated with 
 the arrangement by increase of the minimal elements  of blocks
 are sampling consistent as in Definition \ref{def1}
(this observation is due to Donnelly and Joyce \cite{Donnelly}). 
The last assertion is just another formulation of Theorem \ref{ESForder}.

\subsection{Composition structure derived from the zero set of BM}

Consider the process 
$(B_t, t\geq 0)$ of Brownian motion (BM), and let ${\cal Z}=\{t: B_t=0\}$ be the zero set of BM.
The complement ${\mathbb R}\setminus {\cal Z}$ is the union of the excursion intervals, where the BM is away from zero. 
Define ${\cal R}$ as ${\cal Z}$ restricted to $[0,1]$.
There is a meander 
gap between the last zero of BM and $1$,
caused by an incomplete
excursion   abrupted at $t=1$, but to the left of the meander the set $\cal R$ is
of the Cantor-type, without isolated points. 
Thus the gaps cannot be simply enumerated from left to right, as in the stick-breaking case.
Since the  BM is a recurrent process with the strong Markov property,
the set of zeroes to the right of 
the generic excursion interval
is a shifted distributional  copy of the whole $\cal Z$, 
independent of the part of $\cal Z$ to the left of (and including) the excursion interval.
This implies that $\cal Z$ is a {\it regenerative} set, a property familiar from the elementary renewal theory. 
The scaling property of the BM,
$(c^{-1/2}\,B_{ct})\ed (B_t)$, implies the {\it self-similarity},
 $c{\cal Z}\ed{\cal Z}$ for $c>0$, i.e. the invariance of $\cal Z$ under homotheties.

\par Following Pitman \cite{PitmanBernoulli}, 
consider the composition structure $\varkappa$ 
 derived from  ${\cal R} ={\cal Z}\cap [0,1]$.
To check the deletion property in Definition \ref{def1} it is convenient to modify 
the paintbox model in a way accounting for the self-similarity.

\vskip0.2cm
\noindent
{\bf Modified sampling scheme} Let ${\cal Z}$ be a random self-similar subset of $\mathbb R$.
Fix $n$ and let $X_1<X_2<\ldots$
be the points of a unit Poisson process, independent of $\cal Z$. 
The interval $[0,X_{n+1}]$ is split in components at points of
$\cal Z$,
so we can define a composition $\varkappa_n$ of $n$ by grouping $X_1,\ldots,X_n$ in  clusters within $[0,X_{n+1}]$.  
As $n$ varies, these compositions comprise the same composition structure, as the one induced by the standard 
paintbox construction with ${\cal Z}\cap[0,1]$, because
(i) the vector
$(X_1/X_{n+1},\ldots,X_n/X_{n+1})$ is distributed like the vector of $n$ uniform order statistics $(U_{n:1},\ldots,U_{n:n})$
and (ii) by self-similarity, ${\cal Z}/X_{n+1}\ed\cal Z$.

\vskip0.2cm
Note that, because the locations of `balls' vary with $n$, 
 the model secures a weak consistency of $\varkappa_n$'s, but does not 
produce strongly consistent ordered set partitions ${\rm K}_n$.
Applying the modified scheme in the BM case,
the deletion property is obvious from the regeneration of $\cal Z$ and of the homogeneous Poisson process, these
combined with the self-similarity of  $\cal Z$.


\subsection{Regenerative sets and subordinators}
In the stick-breaking case the regeneration property of the induced composition structure
$\varkappa$ followed from the  observation that  
$\cal R$ remains in a sense the same when its left meander is truncated.  
This could not be applied in the BM case,
since the leftmost gap does not exist.
By a closer look it is seen that 
a weaker property of $\cal R$ would suffice.
 For a given closed ${\cal R}\subset [0,1]$  
define the `droite' point $Z_x:=\min\{{\cal R}\cap[x,1]\},~x\in [0,1]$, which is the right endpoint of the gap covering $x$
(or $x$ itself in the event $x\in {\cal R}$).

\begin{definition}\label{m-reg}
{\rm A random closed set ${\cal R}\subset [0,1]$ is called {\it multiplicatively regenerative}
(m-regenerative for short) if $Z_x$ is independent of  $(1-Z_x)^{-1}(({\cal R}\cap [Z_x,1])  -Z_x )$ and, given $Z_x<1$, 
the distributional identity is fulfilled
$${({\cal R}\cap [Z_x,1])- Z_x\over 1-Z_x}\ed{\cal R}$$
for every $x\in [0,1)$.
}
\end{definition}

\vskip0.2cm
\noindent
{\bf Remark} We do not require explicitly the independence of $[0,Z_x]\cap{\cal R}$ and $(1-Z_x)^{-1}{\cal R}$,
which would correspond to the conventional regeneration property in the additive theory.
In fact, this apparently stronger property follows from the weaker independence property 
due to connection to composition structures. See \cite{RCS} for details
and connection to the bulk-deletion properties of composition structures.
 
\vskip0.2cm

\par  For m-regenerative paintbox $\cal R$ the deletion property of $\varkappa$ 
follows by considering the gap that covers
$U_{n:1}=\min(U_1,\ldots,U_n)$. Then $q(n:\cdot)$ is the distribution of the rank of the largest order statistic in this gap.

\par To relate Definition \ref{m-reg} with the familiar (additive) concept of regenerative set,
 recall that a {\it subordinator} $(S_t,t\geq 0)$ is an increasing right-continuous process with $S_0=0$ and
stationary independent increments (L{\'e}vy process). The fundamental characteristics of subordinator are 
the L{\'e}vy measure $\nut$ on $(0,\infty]$, which controls the intensity and sizes of jumps, 
and the drift coefficient ${\tt d}\geq 0$ responsible for a linear drift component. 
The distribution is determined by means of the Laplace transform
$${\mathbb E}[\exp(-\rho S_t)]=\exp[-t\Phi(\rho)],~~\rho\geq 0,$$
where the Laplace exponent is given by the L{\'e}vy-Khintchine formula
\eq\label{LK}
\Phi(\rho)=\rho{\tt d}+\int_{(0,\infty]} (1-e^{-\rho y})\nut(\dd y).
\en
The L{\'e}vy measure must satisfy the condition $\Phi(1)<\infty$ 
which implies $\nu[y,\infty]<\infty$ and 
also restricts the mass near $0$, 
to avoid  immediate passage of the subordinator to $\infty$.
A positive 
mass at $\infty$ is allowed, in which case  $(S_t)$ (in this case sometimes called 
 killed subordinator) jumps to $\infty$ at some exponential time with rate $\nut\{\infty\}$. Two standard examples of subordinators are 
\begin{enumerate}

\item
Stable subordinators with parameter $0<\alpha<1$, characterised by 
$$\nut(\dd y)={c \alpha\over \Gamma(1-\alpha)}\,y^{-\alpha-1}\dd y\,,~~~{\tt d}=0,~~~\Phi(\rho)=c\rho^\alpha.$$
\item Gamma subordinators with parameter $\theta>0$, characterised by 
$$
\nut(\dd y)=c y^{-1}e^{-\theta x}\dd y, ~~~{\tt d}=0,~~~\Phi(\rho)= c \log(1+\rho/\theta). 
$$
\end{enumerate}
The constant $c>0$ can be always eliminated by a linear time-change.
\par Let $\widetilde{\cal R}=\{S_t,t\geq 0\}^{\rm cl}$ be the closed range of a subordinator. 
By properties of the increments, $\cal R$ is {\it regenerative}: for $Z_y$ the `droite' point at $y>0$,
conditionally given $Z_y<\infty$, the random set $(\widetilde{\cal R}-Z_y)\cap [0,\infty]$ is distributed like $\widetilde{\cal R}$ 
and is independent of $[0,Z_y]\cap\widetilde{\cal R}$ and $Z_y$. 
Also the converse is true: by a  result of Maisonneuve  \cite{Molchanov}
every regenerative set is the closed range of some subordinator, with $(\nut,{\tt d})$ determined uniquely up to a positive multiple.

\par 
Call the increasing process $(1-\exp(-S_t),t\geq 0)$ 
{\it multiplicative} subordinator, and let 
${\cal R}=1-\exp(-\widetilde{\cal R})$ be its range. 
The regeneration property of $\widetilde{\cal R}$ readily implies that $\cal R$ is m-regenerative.
As time passes, the multiplicative subordinator proceeds from $0$ to $1$, thus it is natural to adjust 
the L{\'e}vy measure to the multiplicative framework by transforming $\nut$, by the virtue of  $y\to 1-e^{-y}$, in some 
measure $\nu$ on $(0,1]$, which accounts now for a kind of continuous-time stick-breaking.
We shall still call $\nu$ the  L{\'e}vy measure where there is  no ambiguity.
In these terms the
L{\'e}vy-Khintchine formula becomes
\eq\label{m-LK}
\Phi(\rho)=\rho{\tt d}+\int_0^1 \{1-(1-x)^\rho\}\nu(\dd x).
\en
For integer $1\leq m\leq n$ introduce also the binomial moments of $\nu$
$$
\Phi(n:m)=
{n\choose m} \int_0^1 x^m (1-x)^{n-m}\,{\nu}({\rm \,d}x)+
1(m =1) \,n{\tt d}
$$
(where $1(\cdots)$ stands for indicator), so that $\Phi(n)=\sum_{m=1}^n \Phi(n:m)$.
According to one interpretation of (\ref{LK}), $\Phi(\rho)$ is the probability rate at which 
the subordinator passes through
independent exponential level with mean $1/\rho$. 
Similarly,  $\Phi(n)$ is the rate at which the multiplicative subordinator passes through $U_{n:1}$
 and
$\Phi(n:m)$ is the rate to jump from below $U_{n:1}$ to a value between $U_{n:m}$ and $U_{n:m+1}$.
From this, the probability that the first passage through $U_{n:1}$ covers $m$ out of $n$ uniform points is equal to
\eq\label{dm-Phi}
q(n:m)={\Phi(n:m)\over \Phi(n)}, 
\en
which is the general representation for decrement matrix of a 
regenerative composition structure associated with m-regenerative set. 
The proper case corresponds to the zero drift, ${\tt d}=0$, then passage through a level can 
only occur by a jump.

\par In the case of finite $\nut$ and ${\tt d}=0$ the subordinator is a compound Poisson process with no drift.
Scaling $\nu$ to a probability measure, the range of $(1-\exp(-S_t), t\geq 0)$ 
is a stick-breaking set with the generic factor $W$ distributed according
to $\nu$\,; then (\ref{dm-Phi}) becomes  (\ref{stbr}).

\par The connection between regenerative compositions structures and regenerative 
sets also goes in the opposite direction.

\begin{theorem}\label{ReCoTh}
Every regenerative composition structure can be derived by the paintbox construction
from the range of a multiplicative subordinator, whose
parameters $(\nu, {\tt d})$ are determined uniquely up to a positive multiple.
\end{theorem}

\noindent
\proof
Sampling consistency together with the regeneration imply that the first $n$ rows of the minor $(q(n':\cdot), \,n'\leq n)$ 
are uniquely determined by the last row $q(n:\cdot)$ via formulas
\begin{eqnarray}\label{q-sc1}
q(n':m')&=&{q_0(n':m')\over 1-q_0(n':0)}\,\,,~~~~1\leq m'\leq n',\\
\label{q-sc2}
q_0(n':m')&=&\sum_{m=1}^n q(n:m){{n-m\choose n'-m'}{m\choose m'}\over {n\choose n'}}\,,~~~0\leq m'\leq n'.
\end{eqnarray} 
Think of $\varkappa_n$ 
as allocation of $F_n$ balls in 
box labeled $B$, and $n-F_n$ balls in other boxes. 
Formula (\ref{q-sc2}) gives the distribution of the number of balls remaining in $B$
after $n-n'$ balls have been removed 
at random  without replacement, 
with account of the possibility $m'=0$ that $B$ may become empty.
Formula (\ref{q-sc1}) says that  the distribution of $F_{n'}$ is the same as that of the number of balls which
remain  in $B$ conditionally given that at least one ball remains.
This relation of $F_n$ and $F_{n'}$ is counter-intuitive, because
sampling may eliminate the first block of $\varkappa_n$ completely
(equivalently, the first block of ${\rm K}_n$ may  have no representatives in $[n']$).
\par Invoking the simplest instance of (\ref{q-sc1}), with $n'=n-1$, we have
\begin{eqnarray*}
{q(n:m)\over 1-q(n+1:1)/(n+1)}
={m+1\over n+1}\,q(n+1:m+1)
+{n+1-m\over n+1}\,q(n+1:m).
\end{eqnarray*}
This is a nonlinear recursion, but 
passing to formal homogeneous variables
$\Phi(n:m)$ and using the
substitution   $q(n:m)=\Phi(n:m)/\Phi(n)$ with $\Phi(n):=\sum_{m=1}^n\Phi(n:m)$
results in  
the linear relation
$$\Phi(n:m)={m+1\over n+1}\Phi(n+1:m+1)+{n-m+1\over n+1}\Phi(n+1:m).$$
Equivalently, in terms of the iterated differences
$$\Phi(n:m)={n\choose m}\sum_{j=0}^m (-1)^{j+1}{m\choose j}\Phi(n-m+j), ~~1\leq m\leq n.$$ 
The positivity condition $\Phi(n:m)\geq 0$ implies that the sequence $(\Phi(n), n\geq 0)$
(where $\Phi(0)=0$)
must be
{\it completely alternating} \cite{Berg}, i.e. its iterated differences have alternating signs.
The latter also means
that the difference sequence
$(\Phi(n+1)-\Phi(n),~n\geq 0)$ is completely monotone, hence 
by the famous Hausdorff theorem
$\Phi(n)$'s are representable as 
 moments of some finite measure 
on $[0,1]$.
From this
(\ref{m-LK}) follows for integer values of $\rho$ with some $(\nu,{\tt d})$.
The latter secures (\ref{m-LK}) for arbitrary $\rho>0$ by the uniqueness of interpolation.
\par Interestingly, the  argument 
only exploits a recursion on $q$, hence
avoids explicit limit transition
from $\varkappa$ to $\cal R$, as one could expect by analogy with Theorem \ref{Thm1}.
See \cite{RCS, Dong, Freedman} for variations.
\endpf

\par
We  can also view $F(t)=1-\exp(-S_t)$ as a random distribution function 
on ${\mathbb R}_+$ and to construct a composition by sampling from $F$, as in the species sampling problem.
These {\it neutral to the right priors} have found applications
in  Bayesian statistics 
\cite{James}.

\vskip0.3cm

\noindent
{\bf Additive paintbox}
It is sometimes convenient  to induce regenerative composition structures using a  subordinator $(S_t)$
to create the gaps.
Then, independent unit-rate exponential variables $E_1,E_2,\ldots$ should be 
used in the role of balls, instead of  uniform $U_j$'s in the multiplicative framework.

\par Formula (\ref{dm-Phi}) can be re-derived by appealing to the {\it potential measure} ${\cal u}$ of 
subordinator. 
Heuristically, think of ${\cal u}(\dd y)$ as of probability to visit location $y$
at some time, and of $\nut$ as distribution of size of a generic jump  of the subordinator.
The probability that the first part $m$ of composition is created by visiting $y$ by a jump of given size $z$ is then the product of
$(1-e^{-ny}){\cal u}(\dd y)$ and  $(1-e^{-z})^m e^{-(n-m)z}\nut(\dd z)$. Taking into account the formula for Laplace 
transform of the potential measure \cite{BertoinL}
$$\int_0^\infty (1-e^{-\rho y}){\cal u}(\dd y)={1\over\Phi(\rho)}\,,$$
we arrive at (\ref{dm-Phi}) by integration.
The compensation formula for Poisson processes \cite[p. 76]{BertoinL}
is needed to make this argument rigorous.

\par The advantage of working with ${\mathbb R}_+$ is that the regeneration property  involves no scaling.
A disadvantage is that the asymptotic frequency of balls within the walls $(a,b)$ is the exponential probability $e^{-b}-e^{-a}$,
as compared to the size of gap in the multiplicative representation on $[0,1]$.
\par In particular, for Ewens' composition structures the subordinator  $(S_t)$ is a compound Poisson process with the jump distribution
exponential$(\theta)$, so the range $\widetilde{\cal R}$ of $(S_t)$  is a homogeneous Poisson point process with density $\theta$,
and $\cal R$ is inhomogeneous  Poisson point process with density $\theta/(1-x)$ on $[0,1]$.

\vskip0.2cm
\begin{example}\label{gammageom}
{\rm Consider 
the infinite L{\'e}vy measure on $[0,1]$ with density
$\nu(\dd x)=x^{-1}(1-x)^{\theta-1}\dd x$. Denoting 
$h_\theta(n)=\sum_{k=1}^n {1\over \theta+k-1}$ the generalised harmonic numbers we compute
$$q(n:m)={n! (\theta)_{n-m}\over m(n-m)!(\theta)_nh_\theta(n)},~~~p^\circ(\lambda_1,\ldots,\lambda_k)={f^\circ(\lambda^\circ)\over (\theta)_n}\prod_{j=1}^k
{1\over h_\theta(\Lambda_j)}.$$
This composition structure appears as the limit of stick-breaking compositions structures (\ref{2parbeta}) as $\gamma\to 0$.
Although the 
CPF looks very similar to Ewens' (\ref{ESF}), there is no simple 
product formula for the associated partition structure, even in the case 
$\theta=1$. 
}
\end{example}
 
\begin{example}\label{hook}{\rm  (Regenerative hook compositions)}
{\rm 
Hook composition structures are induced by killed pure-drift subordinators with 
$\nu(\dd x)=\delta_1(\dd x)$ and ${\tt d}\in [0,\infty].$
They have
decrement matrices with the only nonzero entries
$$q(n:n)={1\over 1+n{\tt d}},~~q(n:1)={n{\tt d}\over 1+n{\tt d}}.$$
The compositions $\varkappa_n$ only assume values like
$(1,1,\ldots,1,m)$. 
Ferrer's diagrams of the 
associated partitions
$(m,1,1,\ldots,1)^\darr$,
are $\Gamma$-shaped hooks.
}

\end{example}

\par The hook compositions bridge 
between the pure-singleton composition (with   ${\cal R}=[0,1]$) and 
the trivial one-block composition (with ${\cal R}=\{0,1\}$).
For arbitrary composition structure with some L{\'e}vy exponent $\Phi$ we can construct a similar deformation by adding  atomic component
$\beta\delta_1(\dd x)$ to the L{\'e}vy measure; this results in a family of decrement matrices  
\eq\label{betaki}
{\Phi(n:m)+\beta\, 1(m=n)\over \Phi(n)+\beta},
\en
with the one-block composition appearing in the limit $\beta\to\infty$.

\vskip0.2cm
\noindent
{\bf Sliced splitting} We introduce another kind of parametric deformation of a subordinator. 
Let $S=(S_t)$ be a subordinator with range $\widetilde{\cal R}$, and let 
$X_1<X_2<\ldots$ be the points of homogeneous Poisson process with density $\theta$.
Take $S(j)$,  $j\geq 0$, to be independent copies of $S$, also independent of 
the Poisson process. 
We construct the path of {\it interrupted subordinator} $S^{(\theta)}$ by shifting and glueing  pieces of $S(j)$'s in one path.
\par Run $S(0)$ until the passage through level $X_1$ at some time $T_1$, so
 $S_{T_{1-}}(0)<X_1\leq S_{T_1}(0)$. Leave the path of the process south-west of the point $(T_1,X_1)$ as it is, and cut the rest 
north-east part of the path.
At time  $T_1$ start the process $(S_t(1)+X_1,~t\geq T_1)$ and let it running until passage through $X_2$.
Iterate, creating partial paths running from $(T_j,X_j)$ to $(T_{j+1}, X_{j+1})$.
From the properties of subordinators and Poisson processes,
one sees that $S^{(\theta)}$ is indeed a subordinator.

\par The range of $S^{(\theta)}$,
$$
\widetilde{\cal R}^{(\theta)}:=\bigcup_{j\geq 0}\ [X_j,X_{j+1})\cap (X_j+ \widetilde{\cal R}(j)),
$$ 
can be called sliced splitting.
First ${\mathbb R}_+$ is split at locations $X_j$, then each gap $(X_j,X_{j+1})$ is further split at 
points of $\widetilde{\cal R}(j)\cap (X_j,X_{j+1})$ where $\widetilde{\cal R}(j)\ed \widetilde{\cal R}$ are i.i.d.

\par 
The range of $1-\exp(-S_t^{(\theta)})$ can be constructed 
by a similar fitting in the gaps
between the points $1-\exp(-X_j)$, which are the atoms of a Poisson point process with density $\theta/(1-x)$.

\par Denote, as usual, $\Phi, \nut$ the characteristics of $S$, and 
$\Phi_\theta, \nut_\theta$ the characteristics of $S^{(\theta)}$. 
Then we have
\eq\label{2level}
\Phi_\theta(\rho)={\rho\over \rho+\theta}\Phi(\rho+\theta),~~~~\nut_\theta[y,\infty]=e^{-\theta y}\nut[y,\infty].
\en
To see this, a heuristics is helpful to guess the passage rate through exponential level.
Denote $E_\rho,E_\theta$ independent exponential variables with parameters $\rho,\theta$.
The process $S^{(\theta)}$ passes the level $E_\rho$ 
within infinitesimal time interval $(0,t)$ when $S_t>\min(E_\rho,E_\theta)=E_\rho$.
The inequality $S_t>\min(E_\rho,E_\theta)\ed E_{\rho+\theta}$ occurs with probability $\Phi(\rho+\theta)t+o(t)$, and  
 probability of the event $E_\theta<E_\rho$ is $\rho/(\rho+\theta)$.

\vskip0.2cm
\noindent
{\bf The Green matrix}
For  a sequence of compositions $\varkappa=(\varkappa_n)$ which, in principle, need not be consistent in any sense
we can define $g(n,j)$ as the probability that a `1' stays in position $j$ of the binary code of
$\varkappa_n$. 
That is to say, $g(n,j)$ is the probability that the parts of $\varkappa_n$ satisfy $\lambda_1+\ldots+\lambda_{i-1}=j-1$ for some $i\geq 1$.
Call $(g(n,j), 1\leq j\leq n, n\in \Nat)$ {\it the Green matrix} of $\varkappa$.
For $\varkappa$  a regenerative composition structure, 
$g(n,j)$
is  the probability that the Markov chain $Q_n^\darr$ ever visits state $n+1-j$, and we have an explicit  formula
in terms of the Laplace exponent
(see \cite{RCS})
\eq\label{Green}
g(n,n-j+1)=\Phi(j){n\choose j}\sum_{j=0}^{j-1}{j-1\choose i} {(-1)^i\over \Phi(j+i)}.
\en

\subsection{Regenerative compositions from the two-parameter family}
\label{twoparam}

Let $\pi$ be the two-parameter partition structure with PPF (\ref{2par}).
Sometimes notation PD$(\alpha,\theta)$ is used for
the law of frequencies $\cal S$, where PD stands for Poisson-Dirichlet, and sometimes this law is called
Pitman-Yor prior after \cite{PY97}. 
Formulas for  PD$(\alpha,\theta)$  are difficult, but the sequence of frequencies in size-biased order 
can be obtained by inhomogeneous stick-breaking scheme (\ref{sbreak}) with $W_j\ed {\rm beta}(1-\alpha,\theta+j\alpha)$.

\par We will see that for $0 \leq \alpha < 1 $ and $\theta \geq  0$ and only for these values of the parameters
the parts of $\pi$ can be arranged 
in a regenerative composition structure.

\par Define a (multiplicative) L{\'e}vy measure $\nu$ on $[0,1]$ by  the formula 
for its right tail
\eq
\label{althmeas}
\nu [x,1] = x^{-\alpha} (1-x)^{\theta}.
\en
The density of this measure is a mixture of two beta-type densities, and in the case $\theta=0$ 
there is a unit  atom at $1$.
The associated Laplace exponent is
\eq
\label{althlap}
\Phi (\rho) = \rho {\rm B}(1-\alpha, \rho + \theta) = {\rho \Gamma(1-\alpha) \Gamma(\rho + \theta)
\over \Gamma(\rho + 1 - \alpha + \theta)}\,,
\en
and the binomial moments are
$$
\Phi(n:m) = { n \choose m }\big( { \alpha {\rm B}(m - \alpha, n - m+1 + \theta) + 
\theta {\rm B}( m + 1 - \alpha, n - m + \theta) }\big),
$$
so there exists a regenerative composition structure $\varkappa$ with the decrement matrix  
\begin{equation}\label{q-EP}
q(n:m) ={ \Phi(n:m) \over \Phi(n)}
= 
{ n \choose m } {(1-\alpha)_{m-1} \over ( \theta + n - m)_m} { ( (n - m ) \alpha + m \theta ) \over n }.
\end{equation}
It is a good exercise in algebra to show 
that  the symmetrisation (\ref{symm1}) of the product-form
 CPF  with  decrement matrix (\ref{q-EP}) is indeed the two-parameter PPF (\ref{2par}).


\par
Like their unordered counterparts, the two-parameter regenerative compositions have many 
interesting features. 
Three subfamilies are of special interest and, as the experience shows, should be always analysed first.

\vskip0.2cm
\noindent
{\bf Case $(0,\theta)$ for $\theta \ge 0$.} 
This is the ESF case (\ref{ESF}), with 
$\nu$ being the beta$(1,\theta)$ distribution.
The blocks of composition appear in the size-biased order, the gaps of ${\cal R}^c$ too.

\vskip0.2cm
\noindent
{\bf Case $(\alpha,0)$ for $0 < \alpha <1$.}
\label{alpha-zero}
In this case 
$$\nu({\rm d}x)=   \alpha x^{-\alpha -1} {\rm d}x+ \delta_1 ({\rm d}x)$$
is an infinite measure with a unit atom at $1$.
The  composition structure is  
directed by ${\cal R}= {\cal Z}\cap [0,1]$,
where $\cal Z$ is the range of stable subordinator.
On the other hand, $\cal R$ can be also obtained 
as the range of multiplicative subordinator $(1-\exp(-S_t),t\geq 0)$, where $(S_t)$ 
is the subordinator with L{\'e}vy measure
$$
\nut(\dd y)=\alpha(1-e^{-y})^{-\alpha-1}e^{-y}\dd y+\delta_\infty(\dd y).
$$
The product formula (\ref{produ}) specialises to 
$$
p^\circ(\lambda_1,\ldots,\lambda_k)=f^\circ(\lambda^\circ)\lambda_{k} \alpha^{k-1}\prod_{j=1}^{k}{(1-\alpha)_{\lambda_j-1}\over \lambda_j!}.
$$
This composition structure was introduced in \cite{PitmanBernoulli}, where $\cal Z$ was realised as the 
zero set of a Bessel process of 
 dimension
$2 - 2 \alpha$. For $\alpha = 1/2$ this is the zero set of BM.

\par The decrement matrix $q$ in this case 
has the special property that
there is a probability distribution $h$ on the positive integers such that 
\eq
\label{qnf}
q(n:m) = h(m) \mbox{ if } m < n \mbox{ and } q(n:n) = 1 - \sum_{m = 1}^{n-1}h(m).
\en
This means that $1$'s in 
the binary code of $\varkappa_n$ 
can be identified 
with the set of sites within $1,\ldots,n$ visited by 
a positive random walk on integers (discrete renewal process), with the initial state $1$.
Specifically,
\begin{equation}\label{h-step}
h(m)= { \alpha (1 - \alpha)_{m-1} \over m! },
\end{equation}
and  $q(n:n)= (1 - \alpha)_{n-1} /(n-1)!$.

\par The arrangement of parts of $\pi_n$ in a composition is obtained by
 placing a size-biased part of $\pi_n$ in the {\it last} position in $\varkappa_n$, then by shuffling the remaining parts 
uniformly at random to occupy all other positions.
Exactly the same rule applies on the paintbox level:
for  $\cal S$ following PD$(\alpha,0)$, 
a term is chosen by the size-biased pick and attached to 1 as the meander, 
then the remaining gaps are arranged in the exchangeable order.

\vskip0.2cm
\noindent
{\bf Case $(\alpha,\alpha)$ for $0 < \alpha <1$}. 
The associated
regenerative set has zero drift and the L{\'e}vy measure
$$
\nut({\rm d} y)
= \alpha (1-e^{-y})^{-\alpha-1} e^{-\alpha y}{\rm\, d}y\qquad y\geq 0\,,
$$
this is the zero set of an Ornstein-Uhlenbeck process.
The corresponding range of
multiplicative subordinator 
can be realised as the zero set of a Bessel bridge of
dimension $2 - 2 \alpha$; in the case $\alpha=1/2$ this is the Brownian bridge.

The parts of $\varkappa_n$ are identifiable with 
 the  increments of a random walk with the same step distribution $h$ as in (\ref{h-step}) for the $(\alpha,0)$ case, 
but now
conditioned on visiting the state $n+1$.
The CPF is
\begin{equation}\label{alpha-alpha}
p^\circ_{\alpha,\alpha}(\lambda_1,\ldots,\lambda_k)=f^\circ(\lambda^\circ){\alpha^k \over (\alpha)_n} \prod_{j=1}^k   (1 - \alpha)_{\lambda_i - 1} .
\end{equation}
This function is symmetric for each $k$, which implies that the parts of each $\varkappa_n$ are in the exchangeable random
order.  
This confirms the known fact that the excursion intervals of a Bessel bridge appear in 
exchangeable order.

\par Due to symmetry, the  transition rule from $\varkappa_n$ to $\varkappa_{n+1}$ is a simple variation of the
 Chinese restaurant scheme.
Now the tables are ordered in a row. Given  $\varkappa_n=(\lambda_1,\ldots,\lambda_k)$, 
customer $n+1$ is placed at one of the existing tables with chance $(\lambda_j-\alpha)/(n+\alpha)$ as usual, 
and when a new table is to be occupied, this table is placed with equal probability to the right, to the left or in-between any two 
of $k$ tables occupied so far.

\vskip0.2cm
\par In the case $(\alpha,0)$, there is a right meander appearing in consequence of killing at rate $\nut\{\infty\}=1$.
Removing the atom at $\infty$ yields another $m$-regenerative set (not in the two-parameter family) obtained by (i) splitting 
$[0,1]$ using beta$(1,\theta)$ stick-breaking, (ii) fitting in each gap $(Y_{j-1},Y_j)$ a scaled copy of the $(\alpha,0)$ $m$-regenerative
set. 
The decrement matrix is (\ref{betaki}), with $\Phi$ like for the $(\alpha,0)$ m-regenerative set  and $\beta=-1$.
A dicrete counterpart, $\varkappa_n$, is a path of a random walk with reflection, 
but CPF has no simple formula.

\vskip0.2cm
\noindent
{\bf The m-regenerative  set
 with parameters $0<\alpha<1,\, \theta>0$} is  constructable from the sets $(0,\theta)$ and 
$(\alpha,0)$  by sliced splitting. 
To define a  multiplicative version of the two-level paintbox, to have   a relation like (\ref{2level}),
 first split $[0,1]$ 
at points $Y_i$ of the Poisson process with density $\theta/(1-x)$ as in the Ewens case, 
(recall that this  is 
the same as  stick-breaking with beta$(1,\theta)$ factor $W$). 
Then for each $j$ choose an independent copy of the $\alpha$-stable regenerative set
starting at $Y_{j-1}$ and abrupted at $Y_j$,
and use this copy  to split 
$(Y_{j-1},Y_j)$.

\par The resulting m-regenerative set corresponds to  the 
$(\alpha,\theta)$ composition structure, so (\ref{2level}) becomes
$$
\Phi_{\alpha,\theta}(\rho)={\rho\over \rho+\theta}\Phi_{\alpha,0}(\rho+\theta),
$$
which is trivial to check.
As another  check,
observe that
the structural distribution beta$(1-\alpha,\alpha+\theta)$ is the Mellin convolution of beta$(1,\theta)$ and beta$(1-\alpha,\alpha)$,
as it must be for the two-level splitting scheme.
\par The construction is literally the same on the level of finite compositions. First 
a regenerative Ewens $(0,\theta)$ composition of $n$ is constructed, then each part is independently split in a sequence of parts according 
to the rules of the regenerative $(\alpha,0)$-composition.

\vskip0.2cm
\par The arrangement problem for general $(\alpha,\theta)$ was settled recently in \cite{PitmanWinkel}.
Note that every sequence $r_1,r_2,\ldots$ of {\it initial ranks} $r_j\in [j]$ defines uniquely a
total order on $\Nat$, by  placing $j$ in position $r_j$ relatively to $1,\ldots,j$.
For instance, the initial ranks $1,2,1,3,\ldots$ encode a total order in which the arrangement of set $[4]$ is 
3~1~4~2  (1 is ranked 1 within [1], then 2 is ranked 1 within [2], then 3 is ranked  1 within [3], then 4 is ranked 3 within [4], $\ldots$).
For $\eta\in [0,\infty]$, consider a probability distribution for $(r_1,r_2,\ldots)$ under which $r_j$'s are independent, the probability 
of $r_j=j$ is $\eta/(\eta+j)$ and the probability of $r_j=i$ is $1/(\eta+j)$ for every $i< j$.
Pitman and Winkel \cite{PitmanWinkel} show that to arrange ${\cal S}\ed{\rm PD}(\alpha,\theta)$ in regenerative 
paintbox one should (i) first label the frequencies in the size-biased order, (ii) then, independently, arrange the collection
of frequencies by applying
the arrangement to the lebels, with parameter
$\eta=\theta/\alpha$. For $\alpha=0$, the frequencies will be arranged in the size-biased order
(because for $\eta=\infty$ the relative ranks are $r_j=j$ a.s.); 
for $\alpha=\theta$ this is an exchangeable arrangement of $\cal S$; and for 
$\theta=0$ the arrangement is as for $(\alpha,0)$ partition described above. 
\par The arrangement of blocks of $\pi_n$ in regenerative composition $\varkappa_n$ is analogous, for each $n$.
See \cite{GnPerm} for this and larger classes of distributions on permutations, 
their sufficiency properties and connections to the generalised ESF.

\section{Regenerative partition structures and the problem of arrangement}

We discuss next  connections 
between regenerative composition structures and their associated partition structures.
One important issue is the uniqueness of the correspondence.

\subsection{Structural distributions}

For ${\cal R}^c$ related to $\cal S$ via (\ref{SU}) let $\tilde{P}$ be the size of the gap covering the uniform point
$U_1$,
with the convention that  $\tilde{P}=0$ in the event $U_1\in {\cal R}$. 
We shall understand $\tilde{P}$ as a size-biased pick from $\cal S$, this 
agrees with the (unambiguous) definition in the proper case and extends it when the sum of positive frequencies 
may be less than 1.
Obviously, the particular  choice of $\cal R$ with gap-sizes $\cal S$ is not important. 

\par The law of $\tilde P$ is known as 
the {\it structural distribution} 
of $\cal S$. Most properties of this distribution readily follow from the fact that it is 
a mixture of discrete measures $\sum_j s_j\delta_{s_j}(\dd x)+\left(1-\sum_j s_j\right)\delta_0(\dd x)$.
In particular,
the $(n-1)$st moment of $\tilde{P}$ is the probability  that $\varkappa_n$ is  the trivial one-block composition $(n)$
or, what is the same, that $\pi_n=(n)^\downarrow$:
$$p(n)= {\mathbb E}[\tilde{P}^{n-1}].$$

\par In general, there  can be many  partition structures which share  the same structural distribution, but 
for the 
regenerative composition structures the correspondence is one-to-one.
Indeed,
 we have
$$ p(n)=q(n:n)={\Phi(n:n)\over \Phi(n)}.$$
With some algebra  a recursion for the Laplace exponent follows
$$
\Phi(n)(p(n)+(-1)^n)=\sum_{j=1}^{n-1}(-1)^{j+1}{n\choose j}\Phi(j),
$$
which shows that the moments sequence $(p(n), n\in \Nat)$ determines $(\Phi(n), n\in \Nat)$ uniquely up to a positive multiple,
hence determines the decrement matrix $q$.
Explicit  
expressions of the entries of $q$ through the $p(n)$'s are complicated,
these are some 
rational functions in $p(n)$'s, for instance
$q(3:2)=(2p(2)-3p(3)+p(2)p(3))/( 1-p(2))$.
Because the moments $p(n)$ are determined by the sizes of gaps and not by  their arrangement,
we conclude that

\begin{theorem}\label{un-dis}
Each partition structure corresponds to at most one regenerative composition structure.
Equivalently, for random frequencies $\cal S$ there exists at most one distribution for a m-regenerative set
$\cal R$ with $({\cal R}^c)^\downarrow\ed {\cal S}$.
\end{theorem}

\par In principle, one can determine if some PPF $p$ corresponds to a regenerative CPF $\pi$ by computing $q$ formally 
from the one-block probabilities
$p(n)$'s,  then checking positivity of $q$, and if it is positive
then comparing the symmetrised PPF (\ref{produ}) corresponding to $q$ with $p$.
This method works smoothly in the two-parameter case.
For the $(\alpha,\theta)$ partition structures the structural distribution is
beta$(1-\alpha,\theta+\alpha)$ and
$$
{\mathbb E}[\tilde{P}^{n-1}]= {\theta(1-\alpha)_{n-1}\over (\theta)_{n}}\,,
$$
see Pitman \cite{CSP}. 
Computing formally $q$ from $p(n)$'s  
we arrive at $q$ coinciding with (\ref{q-EP}).
However, a decrement matrix must be nonnegative, which is not the case for some values of the parameters: 
  
\begin{theorem} Every Ewens-Pitman partition structure with parameters in the range $0 \leq\alpha\leq 1,\, \theta\geq 0$ has
a unique arrangement as a regenerative composition structure.  For other values of the parameters 
such arrangement does not exist.
\end{theorem}

\noindent
Actually, it is evident that any partition structure of the `discrete series' with $\alpha<0$ in (\ref{2par-range})
cannot be regenerative just because the number of parts in each $\pi_n$ 
is  bounded by $-\theta/\alpha$.

\subsection{Partition structures invariant under deletion of a part}

Recalling  (\ref{pd}), (\ref{qq}), partition structures inherit a deletion property from the parent regenerative compositions. 
In this section we discuss the  
reverse side of this connection, 
which puts the regeneration property in the new light.
The main idea is that if a partition structure $\pi$ has a part-deletion property, then the iterated deletion
creates order in a way consistent for all $n$, thus canonically associating  with $\pi$ a regenerative composition structure.

\par Let $\pi$ be a partition structure. 
A {\it random part}  of $\pi_n$ is an integer  random variable $P_n$ which satisfies 
$P_n\in \pi_n$.
The joint distribution of $\pi_n$ and $P_n$ is determined by the PPF and some 
 deletion kernel $d(\lambda^\darr,m)$, which specifies the conditional distribution of $P_n$ given partition 
$\pi_n$ 
\eq\label{d-del}
p (\lambda^\downarrow)d(\lambda,m)=\prob(\pi_n=\lambda^\darr, P_n=m),~~~~|\lambda^\darr|=n.
\en
For each $n=1,2,\ldots$ the distribution of $P_n$ is then
\eq\label{regen-p2}
q(n:m)=\prob(P_n=m)=\sum_{\{\lambda^\darr: |\lambda^\darr|=n,\,\, m\in \lambda^\darr\}} d(\lambda^\darr,m)p^\darr(\lambda^\darr),~~
1\leq m\leq n.
\en
The formulas  differ from (\ref{pd}) and (\ref{qq}) in that now they refer to some abstract `random part' $P_n$ of unordered structure.
The requirement that $P_n$ is a part of $\pi_n$ makes 
$$\sum_{{\rm distinct}~m\,\in\, \lambda^\darr} d(\lambda^\darr, m)=1.$$

\begin{definition}
{\rm 
Call a partition structure $\pi=(\pi_n)$ {\it regenerative} if, for each $n$, there exists a joint distribution for $\pi_n$ and its
 random part $P_n$ such 
that for each $1\leq m<n$ conditionally given $P_n=m$ 
the remaining partition $\pi_n\setminus\{m\}$ of $n-m$ has the same  distribution as $\pi_{n-m}$. 
Call $\pi$ regenerative w.r.t. $d$ if the conditional distribution of 
$P_n$ is specified by $d$ as in {\rm (\ref{d-del})}, for each $n$.
Call $\pi$ regenerative w.r.t.  $q$ if $q(n:\cdot)$ is  the law of $P_n$, which means that
\eq\label{regen-p1}
p(\lambda^\darr)d(\lambda^\darr,m)=q(n:m)p(\lambda^\darr\setminus \{m\}), ~~n=1,2,\ldots.
\en
}
\end{definition}

\vskip0.2cm
\noindent
{\bf Example} (Hook partition structures) 
This is a continuation of Example \ref{hook}.
Call $\lambda^\downarrow$ a {\it hook partition} if only $\lambda_1$ may be larger than $1$, for instance
$(4,1,1,1)^\downarrow$. For every deletion kernel with the property 
$$d(\lambda^\darr, 1)=1 ~~{\rm if~~}1\in \lambda^\darr,$$
it can be shown that the only partition structures regenerative w.r.t. such $d$ are those supported by hook partitions, and
they have 
$q(n:n)=1/(1+n{\tt d}),~~~~q(n:1)=n {\tt d}/(1+n{\tt d})$
for some ${\tt d}\in [0,\infty]$.
A regenerative partition structure is of the hook type  if and only if $p((2,2)^\downarrow)=0$.

\vskip0.2cm
\begin{theorem} If 
 a partition structure is regenerative and satisfies $p((2,2)^\darr)>0$ then
$q$ uniquely determines $p$ and $d$, and $p$ uniquely determines $q$ and $d$.
Equivalently, if 
 a regenerative partition structure is not of the  hook type then the corresponding deletion kernel is unique.
\end{theorem}

\subsection{Deletion kernels of the two-parameter family}

For Ewens' composition structures (\ref{ESF}) the deletion kernel is the size-biased pick
$$d_0(\lambda^\darr,m)={k_m m\over n}\,,~~{\rm where~}~k_m=\{j:\lambda_j=m\},~n=|\lambda|.$$
The factor $k_m$ appears since the kernel specifies the chance to choose {\it one} of the parts
of given size $m$, rather than a particular part of size $m$.
The regeneration of Ewens' partition structures under this deletion operation
was observed by Kingman \cite{KingmanEw} and called non-interference,
in a species sampling context.
Kingman also showed that this deletion property is characteristic: if a partition structure is regenerative
w.r.t.  $d_0$, then the PPF is the the ESF (\ref{ESFpart})
with some $\theta\in [0,\infty]$.

\par For the regenerative composition structures of the two-parameter family (with nonnegative $\alpha, \theta$)
the deletion kernel is one of
\eq\label{d-tau}
d_\tau(\lambda^\darr, m):= {k_m\over n}\,{(n-m)\tau+m(1-\tau)\over (1-\tau+(k-1)\tau)},~~~\tau\in [0,1],
\en
where $k=\sum_m k_m$ and $n=|\lambda^\darr|$.
Kingman's characterisation of the ESF is a special case of a more general result (see Gnedin and Pitman \cite{RPS}):
\begin{theorem} Fix $\tau\in [0,1]$.
The only partition structures that are regenerative w.r.t.
 deletion kernel $d_\tau$
are the $(\alpha,\theta)$ partition structures with 
$$0\leq\alpha\leq 1,~ \theta\geq 0~~~{\rm and~~~}\alpha/(\theta+\alpha)=\tau.$$
\end{theorem}

\noindent
Summarising, three  subfamilies are characterised by:
\begin{enumerate}
\item
The kernel $d_0$ is the size-biased choice; only $(0,\theta)$ partition structures are regenerative w.r.t. $d_0$.
\item
The kernel $d_{1/2}$ is a uniform random choice of a part; only $(\alpha,\alpha)$ partition structures are regenerative w.r.t.
$d_{1/2}$ .

\item
The kernel $d_1$ can be called 
{\it cosize-biased deletion}, as each (particular) part $m\in \lambda^\darr$ is selected with
probability proportional to $|\lambda^\darr|-m$;
only $(\alpha,0)$ partitions are regenerative w.r.t. $d_1$.
\end{enumerate}
For general $\tau$, the kernel is intrinsically related to the Pitman-Winkel arrangement
of blocks with $\zeta=\tau^{-1}-1$, see Section 
\ref{twoparam}.

\subsection{Back to regenerative compositions}

The  framework of regenerative partitions suggests to study three objects: the PPF $p$, 
the deletion kernel $d$ and the distribution
of deleted part $q$.
Naively, it might seem that $d$, which tells us how a part is deleted,
 is the right object to start with, like
in  Kingman's characterisation of the ESF
via the size-biased deletion.
However, apart from the deletion kernels $d_\tau$ for the two-parameter family, and kernels  related to hook partitions 
we do not know examples where the approach based on the kernels could be made explicit.
Strangely enough, to understand the regeneration mechanism for partitions, 
one should ignore for a while
the question  {\it how} a part is deleted, and only focus on  
$q$ which tells us {\it what} is deleted.

\par Fix $n$ and let $q(n:\cdot)$ be an arbitrary distribution on $[n]$.
Consider a Markov {\it $q(n:\cdot)$-chain} on the set of partitions of $n$ by which a partition 
$\lambda^\darr$ (thought of as allocation of balls in boxes) is transformed by the rules: 
\begin{itemize}
\item
choose a value of $P_n$ from the distribution $q(n:\cdot)$, 
\item given $P_n=m$ 
sample without replacement $m$ balls 
and discard the boxes becoming empty, 
\item put these $m$ balls in a newly created box. 
\end{itemize}
Similarly, define a Markov $q(n:\cdot)$-chain on compositions $\lambda^\circ$ of $n$ with  
the only difference that the newly created box is placed in the first position. 
Obviously, the $q(n:\cdot)$-chain on compositions projects to the $q(n:\cdot)$-chain on partitions when the order
of boxes is discarded.

\begin{lemma} If {\rm (\ref{regen-p1})} holds for some fixed 
$n$ and distribution $q(n:\cdot)$
then the law of $\pi_n$ is a stationary distribution for the $q(n:\cdot)$-chain on partitions.
\end{lemma}
\proof
The condition (\ref{regen-p1}) may be written as a stochastic fixed-point equation
$$\pi_n\setminus\{P_n\}\ed \widehat{\pi}_{n-P_n},$$
where $(\widehat{\pi}_{n'}, 1\leq n'\leq n)$ is 
a sequence of random partitions, independent of $P_n$, with 
$\widehat{\pi}_n\ed \pi_n$. The lemma follows since then $\widehat{\pi}_{n-P_n}\cup \{P_n\}\ed \pi_n$.
\endpf
\par There is an obvious parallel assertion about a random composition $\varkappa_n$, which 
satisfies 
$$\varkappa_n\setminus \{F_n\}\ed \widehat{\varkappa}_{n-F_n},$$
where $\setminus$ stands for the deletion of the first part $F_n$ with distribution $q(n:\cdot)$.

\begin{lemma} The unique   stationary distribution of the $q(n:\cdot)$-chain on compositions
is the one by which $\varkappa_n$ follows the product formula for $1\leq n'\leq n$ with 
$q(n':\cdot)$ given by {\rm(\ref{q-sc1})}. 
Symmetrisation of the law of $\varkappa_n$  by {\rm (\ref{symm1})} gives the unique 
stationary distribution of the $q(n:\cdot)$-chain on partitions.
\end{lemma}

\par It follows that if (\ref{regen-p1}) holds for some $n$ then it holds for all $n'\leq n$, with all
$p(\lambda^\downarrow), d(\lambda^\darr,\cdot)$ for $|\lambda^\downarrow|=n'$  uniquely determined by $q(n:\cdot)$ via sampling 
consistency. 
Thus, in principle, 
for partitions of $n'\leq n$ the regeneration property is uniquely determined by
{\it arbitrary} discrete distribution 
$q(n:\cdot)$ through the following steps: find 
first $(q(n':\cdot), n'\leq n)$  from sampling consistency (\ref{q-sc1}), then use
the product formula for compositions (\ref{produ}), then the symmetrisation (\ref{symm1}). With all this at hand,
the deletion kernel can be determined from (\ref{d-del}). 
Letting  $n$ vary, the sampling consistency of all $q(n:\cdot)$'s implies that $q$ is a decrement matrix
of a regenerative composition structure.

\par Starting with $\pi_n$, the deletion kernel determines a Markov chain on subpartitions of $\pi_n$. 
A part $P_n$ is chosen according to the kernel $d$ and deleted, from the remaining partition $\pi_n\setminus\{P_n\}$
another part is chosen according to $d$ etc.
This brings the parts of $\pi_n$ in the {\it deletion order}.

\begin{theorem} Suppose a partition structure $\pi=(\pi_n)$ is regenerative w.r.t. $q$, then
\begin{itemize}
\item[{\rm (i)}] $q$ is a decrement matrix of some regenerative composition structure $\varkappa$,
\item[{\rm (ii)}] $\pi$ is the symmetrisation of $\varkappa$,
\item[{\rm (iii)}]  
$\varkappa$ is obtained from $\pi$ by 
arranging, for each $n$, the parts of $\pi_n$ in the deletion order.
\end{itemize}
\end{theorem}

\par
Thus the regeneration concepts for partition  and composition structures coincide.
It is not clear, however, how to formulate the regeneration property in terms of the unordered frequencies $\cal S$.
The only obvious way is to compute PPF and then check if the PPF corresponds to a regenerative 
CPF. 
Moreover, the deletion kernel may have no well-defined continuous analogue.
For instance, in the $(\alpha, \alpha)$ case $d_{1/2}$ is a uniform random choice of a part from $\pi_n$, 
but what is a `random choice of a term' from the infinite random series $\cal S$ under PD$(\alpha,\alpha)$?

\subsection{More on $(\alpha, \alpha)$ compositions: reversibility}

We have seen that 
the $(\alpha, \alpha)$ composition structures are the only regenerative compositions which
have  parts in the exchangeable order. 
We show now that these structures can be characterised by some weaker properties of reversibility.

\par Every  composition structure 
${\varkappa} $ has a {\it dual\,} 
$ \widehat{\varkappa} $, where each
$\widehat{{\varkappa}}_n $ is the sequence of parts of
${\varkappa}_n $ read in the right-to-left order. For example, the value
$(3,2)$ of $\varkappa_5$ 
corresponds to the value $(2,3)$
of $\widehat{\varkappa}_5$.
If
$ {\varkappa}$ is derived from 
${\cal R}$, then $\widehat{{\varkappa}} $ is derived from  the  reflected paintbox $1 - {\cal R}$.
If both ${\varkappa} $ and $\widehat{{\varkappa}}$ are
regenerative then by the uniqueness (Theorem \ref{un-dis}) they must have the same distribution.
If ${\varkappa}$ is {\it reversible\,}, 
i.e. $\varkappa\ed\widehat{\varkappa}$,
then the first part of $\varkappa_n$ must have the same distribution as its last part.

\begin{theorem}
\label{symm}
Let ${\varkappa}$ be a 
regenerative composition structure. Let  $F_n$ denote the first and $L_n$ the last part of $\varkappa_n$.  
The following conditions are equivalent:
\begin{itemize}
\item[{\rm (i)}] $\prob (F_n = 1) = \prob ( L_n = 1)$ for all $n$;
\item[{\rm (ii)}] $F_n \ed L_n$ for all $n$;
\item[{\rm (iii)}] 
${\varkappa}_n\ed\widehat{\varkappa}_n$ for all $n$ (reversibility), 
\item[{\rm (v)}]
  ${\varkappa}$
is an $(\alpha, \alpha)$-composition structure with some $0\leq \alpha\leq 1$.
\end{itemize}
\end{theorem}
\noindent
\proof Some manipulations with finite differences yield
$$
\prob(F_n=1)=q(n:1)={\Phi(n)-\Phi(n-1)\over \Phi(n)/n}\,\,,~~~
\prob (L_n = 1) = 
n \left[1 - \sum_{ k = 2 }^{n } { n - 1 \choose k - 1 } { (-1)^k \over \Phi (k) } \right].
$$
Equating these probabilities, one arrives at 
$\Phi(n)={(1+\alpha)_{n-1}/(n-1)!}$
where $\Phi(2):=1+\alpha$ and the normalisation $\Phi(1)=1$ is assumed.
The latter is the Laplace exponent corresponding to the $(\alpha, \alpha)$ composition.
\endpf
\vskip0.2cm
\noindent
Invoking the paintbox correspondence, the result implies
\begin{corollary}
For a random closed subset ${\cal R}$ of $[0,1]$, the following
two conditions are equivalent:
\begin{itemize}
\item
[{\rm (i)}] ${\cal R}$ is m-regenerative and 
${\cal R} \ed 1 - {\cal R}$.
\item
[{\rm (ii)}] ${\cal R}$ is distributed like the zero set of a 
 Bessel bridge of dimension $2 - 2 \alpha$, for some $0\leq \alpha\leq 1$.
\end{itemize}
\end{corollary}
\noindent
The degenerate boundary cases with $\alpha=0$ or $1$ are  defined by continuity.

\section{Self-similarity and stationarity}

Self-similarity 
of a random closed set ${\cal Z}\subset{\mathbb R}_+$ is
the condition
$c{\cal Z}\ed{\cal Z},$ ~$c>0$.
The property is a multiplicative analogue of the stationarity property 
(translation invariance) of a random subset of $\mathbb R$, 
as familiar from the elementary renewal theory (see \cite{Molchanov} for a general account).
We encountered self-similarity  in connection with paintboxes 
for $(\alpha,0)$ compositions.  
\par Regenerative 
$(0,\theta)$ compositions can be also embedded in the self-similar framework by passing to duals.
The mirrored  paintbox for the dual Ewens' composition structure
is the stick-breaking set ${\cal R}=\{V_1\cdots V_i,~ i=0,1,\ldots\}$ with i.i.d.  $V_i\ed{\rm beta}(\theta,1)$.
This set is the restriction to $[0,1]$ of a self-similar Poisson point process with density $\theta/y,~y>0$.

\par Introduce the operation of {\it right reduction} as cutting the last symbol of the binary code of  composition.
For instance, the right reduction maps $100110$ to $10011$.

\begin{definition}
{\rm 
A  sequence of random compositions $\varkappa=(\varkappa_n)$ is called {\it right-consistent} if the right reduction
maps $\varkappa_{n+1}$ in a stochastic copy of $\varkappa_n$.
If  $\varkappa$ is a composition structure, we call it
{\it self-similar} if it is right-consistent.
}
\end{definition}

If a sequence of compositions $\varkappa=(\varkappa_n)$ is right-consistent, it can be realised
on the same probability space as a {\it single infinite} random binary string $\eta_1,\eta_2,\ldots$, 
with $\varkappa_n$ being the composition encoded in the first $n$ digits   $\eta_1, \ldots,\eta_n$.
For right-consistent $\varkappa$ the Green matrix is of the form 
$$
g(n,j)=\prob(\eta_j=1), ~~~~1\leq j\leq n,\,\,n=1,2,\ldots
$$
and we shall simply write $g(j)$.

\begin{theorem} A composition structure $\varkappa$ is self-similar iff the paintbox 
$\cal R$ is the restriction to $[0,1]$ of a selfsimilar set $\cal Z$.
In this case $\varkappa$ can be encoded in an infinite binary string.
\end{theorem}
\proof
The `if' part 
is easily shown using 
the modified sampling scheme, as in the BM example.
The `only if' part exploits convergence of random sets as in Theorem \ref{Thm1}.
\endpf

\par 
Arbitrary infinite binary  string  $\eta_1,\eta_2,\ldots$ (starting from $1$) need not correspond to a 
composition structure, because care 
of the sampling consistency should 
be taken. Let us review the $(0,\theta)$ and $(\alpha,0)$ compositions from this standpoint.

\vskip0.2cm
\noindent
{\bf Example.}
{\rm For $\theta>0$ let $\eta_1,\eta_2,\ldots$ be a Bernoulli string with independent digits and
$$
g(j)=\prob(\eta_j=1)={\theta\over j+\theta-1}.
$$
This encodes the dual Ewens' composition structure, with the last-part deletion property.
In the modified sampling scheme,
the role of balls is taken by
a homogeneous Poisson point process, and the  boxes are created by points of an independent self-similar Poisson process.

\par The family of composition structures can be included in a Markov process with $\theta\geq 0$ considered as a continuous
time parameter \cite{Ewens}. 
On the level of paintboxes the dynamics amounts to intensifying Poisson processes, so that within time $\dd\theta$
the Poisson process
${\cal Z}={\cal Z}_\theta$  is superimposed with another independent Poisson process with density $\theta/x$.
This is an instance of sliced splitting, so 
 (\ref{2level}) is in force.
From this viewpoint a special feature is that the $\theta$-splitting are consistently 
defined, also in terms of interrupted subordinators, which are here compound Poisson processes with exponential jumps.

\par Remarkably, the splitting process remains Markovian in terms of the binary codes, and has the dynamics in which every `0' eventually turns in `1' by the rule:
at time $\theta$,  a `0' in the generic position $j$ of the code 
is switching at rate $1/(\theta+j-1)$ to a `1', independently of digits in all other positions.

\vskip0.2cm
\noindent
{\bf Example.}
{\rm For $\alpha\in (0,1)$ let 
$(T_k)$ be a discrete 
renewal process with $T_0=1$ and independent increments with distribution
$$\prob(T_{k+1}-T_k=m)=(-1)^{m-1}{\alpha\choose m}$$
(the case $\alpha=1/2$ is related to the recurrence time of a standard random walk).
For
$\eta_j=1(\cap_{k\geq 0}  \{T_k=j\})$ 
 the sequence $\eta_1,\eta_2,\ldots$ encodes the regenerative $(\alpha,0)$ composition structure.
The Green matrix is 
$g(j)=(\alpha)_{n-j}/(n-j)!.$

\vskip0.2cm

It is known \cite{Young} that no other Bernoulli or renewal strings are sampling consistent, i.e. 
produce composition structures. We shall turn to a larger class of strings with a Markov property, but first 
review a few general features of the self-similar compositions.

\par Let $\tilde{P}_n$ be the size-biased pick from $\varkappa_n$, and $L_n$ be the last part of the composition.
Similarly, let $\tilde{P}$ be the size-biased  gap-length of $\cal R$, and $L$ be the size of the meander gap
adjacent to $1$.

\begin{theorem}\label{LP} 
Let $\varkappa$ be a self-similar composition structure, thus derived from some self-similar 
set $\cal Z$. Then
\begin{itemize}
\item[\rm (i)] $\tilde{P}\ed L$, ~~{\rm and~~(ii)}~~  $\tilde{P_n}\ed L_n$,
\end{itemize}
and the Green matrix is
 $g(j)=\ex(1-\tilde{P})^{j-1}.$

\end{theorem}
\proof
Since reducing the last box by one ball has the same effect as reducing the box chosen by the size-biased pick, the sizes of the boxes 
must have the same distribution. This yields (ii), and (i) follows as $n\to\infty$.
Alternatively, inspecting the gap covering $U_{n:n}$  it is seen that 
$\ex [L^{n-1}]=p^\circ(n)$, the probability of one-block composition, so the moments of $\tilde{P}$ and $L$ coincide.
Similarly,
$\eta_j=1$ in the event $U_{n:1}>\max({\cal Z}\cap [0,U_{n:j}])$.
\endpf  
The identity (ii)
together with a generalisation of a result by Pitman and Yor \cite{PitmanYorRDD} yields a characterisation of structural distributions,
and shows that $\tilde{P}$ has a decreasing density on $(0,1]$.

\begin{theorem} {\rm \cite{selfsim}}
The structural distribution for self-similar composition structure
is of the form
\eq\label{str-di}
\prob(\tilde{P}\in \dd x)={\nu[x,1]\over ({\tt d}+{\tt m})(1-x)}\,\,\dd x      +{{\tt d}\over{\tt d}+{\tt m}}\delta_0(\dd x),~~~~~x\in [0,1],
\en
where ${\tt d}\geq 0$ and $\nu$ is a measure on $(0,1]$ with
$$
{\tt m}:=\int_0^1 |\log(1-x)|\nu(\dd x)<\infty.
$$
\end{theorem}
\noindent
There is no   atom at $0$ iff ${\tt d}=0$ iff
$\cal Z$ has Lebesgue measure zero.

\subsection{Markovian composition structures}

For a time being we switch  to regeneration in the right-to-left order of parts, starting from the last part,
like for the dual Ewens' composition.
This is more convenient in the self-similar context since $0$ is the center of homothety. 
We first modify the deletion property of compositions by allowing a special 
distribution for the first deleted part (which is now the last part of the composition).

\begin{definition}{\rm 
A composition structure is called {\it Markovian} if the CPF is of the product form

\eq\label{produM}
p^\circ(\lambda^\circ)=q^{(0)}(n:\lambda_k)\prod_{j=1}^{k-1} q(\Lambda_j:\lambda_j),~~~~~\Lambda_j=\lambda_1+\ldots+\lambda_j.
\en
where $q^{(0)}$ and  $q$ are two decrement matrices.
}
\end{definition}
\noindent
Similarly to (\ref{produ}), formula (\ref{produM}) says that 1's in the binary  code of $\varkappa$ 
appear at sites $Q_n^\darr(t)+1$ visited by a decreasing Markov chain, with the only new feature  that the
the distribution of the first decrement is determined by
 $q^{(0)}$, and not  by $q$.

\par The counterpart  of Theorem \ref{ReCoTh} for (\ref{produM}) is straightforward.
For $(S_t)$  a subordinator, consider
the process
$(V\cdot \exp(-S_t),~ t\geq 0)$,  where $V$ takes values in $(0,1)$ and is 
independent of $(S_t)$. The range of this process is a m-regenerative set (now with right-to-left regeneration)
scaled by the random factor $V$. Taking this set for paintbox $\cal R$, thus with the meander gap $[V,1]$,
 a Markovian composition structure is induced with
$q(n:m)=\Phi(n:m)/\Phi(n)$ as in (\ref{dm-Phi}), and 
$$ 
q^{(0)}(n:m)=\Phi^{(0)}(n:0)q(n:m)+\Phi^{(0)}(n:m),~~~~~\Phi^{(0)}(n:m):={n\choose m}\ex\{V^{n-m} (1-V)^{m}\}.
$$
Every Markovian composition structure is of this form.

\subsection{Self-similar Markov composition structures}

Let $Q^\uparrow=(Q^\uparrow(t),~t=0,1,\ldots)$ be a time-homogeneous increasing Markov chain on $\Nat$ with $Q^\up(0)=1$.
An infinite  string $\eta_1,\eta_2,\ldots$ is defined as the sequence of sites visited by $Q^\uparrow$ 
$$\eta_j=1(Q^\uparrow(t)=j~~~{\rm for~some~}t).$$
If the string determines some composition structure $\varkappa$, then $\varkappa$
is self-similar.
A composition structure 
is called {\it self-similar Markov} if it has such a binary representation generated by an increasing Markov
chain. 

\par A {\it stationary regenerative set} 
(or stationary Markov \cite{Molchanov})
is the range of a process $(X+S_t, ~t\geq 0)$ where 
$(S_t)$ is a finite mean-subordinator, with L{\'e}vy measure satisfying
$$
{\tt m}=\int_0^\infty y\nut(\dd  y)<\infty,
$$
 drift ${\tt d}\geq 0$  and   the initial value $X$ 
whose distribution is
$$
\prob(X\in \dd y)={\nut[y,\infty]\over {\tt d}+{\tt m}}\,\dd y+{{\tt d}\over {\tt d}+{\tt m}}\delta_0(\dd y)
$$
(unlike $\nu$ in (\ref{str-di}) $\nut$ lives on $(0,\infty)$).

\begin{theorem}{\rm \cite{selfsim}}
A composition structure $\varkappa$ is self-similar Markov if and only if 
${\cal R}=\exp(-\widetilde{\cal R})$, where $\widetilde{\cal R}$ is a stationary regenerative set.
\end{theorem}
The distribution of size-biased pick is then (\ref{str-di}) with $\nu$ the image of $\nut$ under $y\to 1-e^{-y}$.
The Green matrix can be written in terms of the Laplace exponent 
$$
g(j)={1\over {\tt d}+{\tt m}} {\Phi(j-1)\over j-1},~~~{~\rm for~}j>1,~ g(1)=1.
$$
The relation beween this and 
(\ref{Green}) is that the RHS of (\ref{Green}) converges to
$g(j)$ as $n\to\infty$. This fact is analogous to the elementary renewal theorem.


\par 
Like in the regenerative case, the decrement matrices are determined, in principle,
 by the probabilities $(p(n),n\geq 0)$, which are moments of the structural distribution,
whence the analogue of Theorem \ref{un-dis}:

\begin{theorem}
If a partition structure admits arrangement as a self-similar Markov composition structure, then such arrangement is unique in distribution.
\end{theorem}
\noindent

\vskip0.2cm\noindent
{\bf Application to the two-parameter family}
For $0\leq \alpha<1$ and $\theta>0$ let ${\cal R_{\alpha,\theta}}$ 
be the m-regenerative set associated with  $(\alpha,\theta)$ regenerative composition structure,
and let $V$ be an independent variable whose  distribution is beta$(\theta+\alpha, 1-\alpha)$. 
Then the scaled reflected set 
$V\cdot(1-{\cal R}_{\alpha,\theta})$  is associated with a self-similar 
Markov composition structure corresponding to $(\alpha,\theta-\alpha)$ partition struture.
This follows from the  
  stick-breaking representation of the frequencies in size-biased order, with
 independent factors beta$(\theta+ j\alpha, 1-\alpha),~ j=1,2,\ldots$.
The Green function $g$ and transition probabilities for $Q^\uparrow$ can be readily computed.

\par A `stationary' version of the regenerative $(\alpha,\theta)$ composition is the self-similar Markov arrangement of the 
$(\alpha, \theta-\alpha)$ partition. The structural distribution is beta$(1-\alpha, \theta+\alpha)$, which is
also the law of the meander size $1-V$.
Note that $\theta-\alpha$ may assume negative values, hence every partition with $\theta>-\alpha$ has a self-similar Markov
arrangement.
This `rehabilitates' $(\alpha, \theta)$ partitions with $-\alpha<\theta<0$ that lack regeneration literally,
the property appears in a modified form, as stationary regeneration.
If $\theta\geq 0$  then both types of regeneration are valid\footnote{
For `discrete series' of the parameter values, with $\alpha<0$, no  regeneration property can exist, simply because the paintbox 
has uniformly bounded cardinality.}.
\par The $(\alpha, 0)$ composition with left-to-right regeneration is also
self-similar Markov, i.e. has the `stationary' right-to-left regeneration property.
This combination of regeneration properties is characteristic for this class.

\par For the $(\alpha,\alpha)$ partition structure there exists a regenerative arrangement associated with
Bessel bridge, and there is another self-similar Markov arrangement.
The latter is the self-similar version of the regenerative $(\alpha,2\alpha)$ composition.

\par The arrangement of $(\alpha, \theta)$ partition in a self-similar Markov composition structure is the same 
on both paintbox and finite$-n$ level. The size-biased pick is placed at the end, then the rest parts 
are arranged to the left of it as for the dual $(\alpha, \theta+\alpha)$ regenerative structure,
see Section \ref{twoparam}.
Property (i) in Theorem \ref{LP}
holds in the strong sense:
conditionally given the unordered frequencies $\cal S$, 
the length of the meander  is a size-biased pick
(see \cite{PitmanYorRDD}).

\section{Asymptotics of the block counts}

For $\varkappa=(\varkappa_n)$ a regenerative composition structure, let $K_n$ be the number of parts in $\varkappa_n$ 
and let $K_{n,r}$ be the number of parts equal $r$, so that $\sum_r rK_{n,r}=n,~\sum_r K_r=K_n$. 
For instance, in the event $\varkappa_{10}=(2,4,2,1,1)$ we have $K_{10}=5, K_{10,1}=2,  K_{10,2}=2, K_{10,3}=0$ etc. 
The full vector $(K_{n,1},\ldots,K_{n,n})$ is one of the ways to record the partition associated with $\varkappa_n$.
In the species sampling context, $K_n$ is the number of distinct species represented in a sample, hence 
it is often considered as a measure of   diversity.

\par We are interested in the large-$n$ asymptotics of $K_n$ and $K_{n,r}$ for $r=1,2,\ldots$.
This can be called the {\it small-blocks problem}.
Typically the composition will have a relatively few number of large parts of size of order $n$ and many parts of size $r\ll n$, the latter making the 
principal contribution to $K_n$.


\par Unless indicated otherwise, we assume that ${\tt d}=0$ (proper case, no drift) and that $\nut\{\infty\}=0$ 
(no killing, no right meander). 
Then the order of growth of $K_n$ is sublinear, $K_n\ll n$, and $K_n\up\infty$ almost surely. 

\par One general tool is the structural distribution $\sigma$ of the size-biased pick $\tilde{P}$, 
which can be used to compute the expectations via
$$
\ex[K_n]= \int_0^1{1-(1-x)^n\over x}\,\sigma(\dd x),~~~~\ex[K_{n,r}]= {n-1\choose r-1}\int_0^1 x^{r-1}(1-x)^{n-r}\sigma(\dd x).
$$
It is clear from these formulas that the asymptotics of the moments are determined by the behaviour
of $\sigma$ near $0$, because $(1-x)^n$ decays exponentially fast 
on any interval $[\epsilon,1]$.

\par The block counts $K_n, K_{n,r}$ depend only on the partition, and not on the order of the parts.
Nevertheless, the Markovian character of regenerative compositions and the connection with subordinators 
can be efficiently exploited to
study these functionals by methods of the renewal theory. This may be compared with
other classes of partitions studied with the help of local limit theorems:
partitions obtained by conditioning random sums of independent integer variables \cite{ABT},
and partitions derived from conditioned subordinators \cite{PK}.

\par For Ewens' partitions it is well known that $K_n$ is asymptotically normal, with both mean and variance of the 
order of $\log n$ (see \cite{ABT, CSP}). 
In contrast to that, for $(\alpha,\theta)$ partitions with $\alpha>0$ 
the right scale for $K_n$ is $n^\alpha$ ($\alpha$-diversity \cite{CSP}).
These known facts will be embedded in a much more general consideration.

\par The number of parts satisfies a  distributional fixed-point equation
$$
K_n\ed 1+K_{n-F_n}'
$$
where $K_{m}', m\leq n-1,$ are independent of the first part $F_n$ with distribution $q(n:\cdot)$,
and satisfy $K_m'\ed K_m$.
Known asymptotics 
(e.g. \cite{Nein}, \cite{Cutsem}) 
derived from such identities  
do not cover the full range of possibilities 
and require very restrictive moment conditions 
which are not easy to provide
(see however  \cite{GPYII} for one application of this approach).
In what follows we report on the asymptotics which were obtained by different methods, based 
on the connection with subordinators,
poissonisation, methods of the renewal theory,
 and Mellin transform \cite{GPYI, GPYII, BG, GnIksM}.

\par We assume as before the paintbox construction with 
balls $U_{1},\ldots,U_{n}$ and
${\cal R}$ the closed range of a multiplicative subordinator $(1-\exp(-S_t), t\geq 0)$.
In these terms, $K_{n,r}$ is the number of gaps in the range hit by exactly $r$ out of $n$ uniform points, and
$K_n$ is the total number of nonempty gaps.

\vskip0.2cm
\noindent
{\bf Remark}  If the subordinator has positive drift ${\tt d}>0$, then
$K_n\sim K_{n,1}\sim n\,{\rm meas}({\cal R})$ a.s., so singletons make a leading contribution to $K_n$.
The Lebesgue measure of $\cal R$ 
is a random variable  proportional to the exponential functional of the subordinator,
$$
{\rm meas}({\cal R})={\tt d}\int_0^\infty \exp(-S_t)\dd t\,.
$$
\vskip0.2cm
\noindent

\par It is informative to consider the number of parts $K_n$ as the   terminal value of the increasing process $\KK_n:=(\KK_n(t), t\geq 0)$,
where ${\cal K}_n(t)$ 
is the number of parts 
of the subcomposition derived from the configuration of  uniform points not exceeding
 $1-\exp(-S_t)$, i.e. produced by the subordinator within the time $[0,t]$. 
The number of $r$-parts $K_{n,r}$ is the terminal value of another process $\KK_{n,r}:=(\KK_{n,r}(t), t\geq 0)$ which counts $r$-parts, 
but this process is not monotone. 
\par We can think of the subordinator representation of a regenerative composition structure as a coagulation process in which, 
if at time $t$ there are 
$n'$ particles, every $m$-tuple of them is merging to form a single particle at rate $\Phi(n':m)$. The particle emerging from the 
coalescence is immediately frozen\footnote{The dynamics is analogous to that of Pitman-Sagitov $\Lambda$-coalescents, with the difference that 
in the $\Lambda$-coalescents the mergers remain active and keep  coagulating with other existing particles
\cite{lambda}. 
For a class of $\Lambda$-coalescents
 a coupling with compositions was used to explore asymptotics of the coalescent processes  \cite{GnIksM}.}.  
Starting with $n$ particles, ${\cal K}_n(t)$ counts the number of frozen particles at time $t$.

\par The asymptotics in the small-block problem largely depend  on the behaviour of the right tail of the L{\'e}vy measure near $0$. 
If $\nut$ is finite, then simply $\nut[y,\infty]\to \nut[0,\infty]$ as $y\to 0$,  
but if $\nut$ is infinite
it seems difficult if at all possible to make any conclusions without the following  assumption.  
 
\vskip0.2cm
\noindent
{\bf Assumption of regular variation}
We shall suppose  that $\nut$ satisfies the condition of regular variation
\eq\label{regvar}
\nut[y,\infty]\sim\ell(1/y) y^{-\alpha}~~~~y\darr 0,
\en
where the index satisfies $0\leq \alpha\leq 1$ and $\ell$ is a function of slow variation at $\infty$,
i.e. $\ell$ satisfies $\ell(t/y)/\ell(1/y)\to 1$ as $y\to 0$ for all $t>0$.
 
\vskip0.2cm

\par Note that the assumption is satisfied in the case of finite $\nut$. 
By
the monotone density version of Karamata's  Tauberian theorem \cite{BGT}, 
for $0\leq \alpha<1$ the condition (\ref{regvar}) is equivalent to the asymptotics of the Laplace exponent
$$\Phi(\rho)\sim \Gamma(1-\alpha)\rho^\alpha\ell(\rho),~~~~~\rho\to\infty.$$ 
Qualitatively different asymptotics are possible.
Very roughly, 
the whole spectrum can be divided in  
the following cases, each requiring separate analysis.

\begin{itemize}

\item The finite L{\'e}vy measure case.
This is the case of stick-breaking compositions, with $(S_t)$ a compound Poisson process.  

\item The slow variation  case with $\alpha=0$ and $\ell(y)\to \infty$ as $y\to\infty$.
Typical example: regenerative compositions associated with gamma subordinators.

\item The proper regular variation case with $0<\alpha\leq 1$.
Typical example: composition associated with  $\alpha$-stable subordinator (with $0<\alpha<1$).
\end{itemize}

\par One principal difference between the cases of (proper) regular and slow variation
is in the time scales at which major growth and variability of $\KK_n$ occur. 
In the case $\alpha>0$ all $\KK_n(t), \KK_{n,r}(t)$ are of the same order as $K_n$, whereas in the case $\alpha=0$ 
we have
$\KK_n(t)\ll K_n$.

\subsection{Stick-breaking compositions}

In the case of finite L{\'e}vy measure we scale  $\nut$  to  a probability measure.
Then $\nut$ is the distribution of $-\log(1-W)$, where $W$ is the generic stick-breaking factor.
 Introduce the moments 
\eq\label{momstbr}
{\tt m}:=\ex[-\log(1-W)],~~~\sigma^2:={\rm Var}[\log(1-W)],~~~{\tt m}_1:=\ex[-\log W],
\en
which may be finite or infinite.
\par 
Let $M_n$ be the index of the rightmost occupied gap, which contains
the maximum order statistic $U_{n:n}$. 
Roughly speaking, stick-breaking implies a fast exponential decay of the sizes of gaps, 
hence
one can anticipate a cutoff phenomenon: empty gaps can occur only in a range close to $U_{n:n}$.
From the extreme-value theory we know that $-\log(1-M_n)-\log n$ has a limit distribution of the Gumbel type, 
thus $M_n$ can be approximated by the number of jumps of $(S_t)$ before crossing level $\log n$.

\par It should be noted that exponential decay of nonrandom frequencies, like for the geometric distribution, 
implies oscillatory asymptotics in the occupancy problem \cite{BGY}, \cite{BG2}.
 By stick-breaking the oscillations do not appear since
the main variability comes due to randomness
in frequencies themselve, so the variability coming from sampling is dominated.

\par Consider a renewal process 
with distribution for spacings like that of $-\log(1-W)$.
If the moments are finite, ${\tt m}<\infty, ~\sigma^2<\infty$,  then
a standard result from the renewal theory implies that
the number of renewals on $[0,\log n]$ is approximately normal for large $n$, with the expected value asymptotic to
$(\log n)/{\tt m}$. The same is valid for $M_n$, and under the additional
assumption ${\tt m}_1<\infty$ also for  $K_n$ 
(see \cite{BS}).
Under weaker assumptions on the moments,  the
possible asymptotics 
correspond to other limit theorems of  
renewal theory,  as shown in  \cite{GnIksM}:

 \begin{theorem}\label{main1} Suppose the distribution of $-\log(1-W)$ is nonlattice with
${\tt m}_1<\infty$.
The following assertions are equivalent.
\begin{enumerate}
\item[\rm (i)] There exist constants $a_n, b_n$
with $a_n>0$ and $b_n\in{\mathbb R}$ such that, as $n\to\infty$, the variable
$(K_n-b_n)/a_n$ converges weakly to some
non-degenerate and proper distribution.
\item[\rm (ii)] The distribution $\nut$ of $-\log(1-W)$
either belongs to the domain of attraction of a stable law, or the function
$\nut[x,\infty]$ slowly varies as $x\to\infty$.
\end{enumerate}
Furthermore, this limiting distribution of $(K_n-b_n)/a_n$ is as follows.
\begin{enumerate}
\item[\rm (a)]
 If $\sigma^2<\infty$,
then for
$b_n={\tt m}^{-1}\log n$ and $a_n=({\tt m}^{-3}\sigma^2\log n)^{1/2}$
the limiting distribution is standard normal.
\item[\rm (b)]
If
$\sigma^2=\infty$ and
$$\int_x^1 (\log y)^2\, 
\prob(1-W\in {\rm d}y)
\sim \ell(-\log x) \ \ \text{as} \
x\to 0,$$ for some $\ell$ slowly varying at $\infty$, then for
$b_n={\tt m}^{-1}\log n$, $a_n={\tt m}^{-3/2}c_{\lfloor\log n\rfloor}$ and
$c_n$  any sequence satisfying
$\lim_{n\to\infty}n\ell(c_n)/c_n^2=1$,
the limiting distribution is standard
normal.
\item[\rm (c)]
Assume that the relation
\begin{equation}\label{domain1}
\prob(1-W\leq x)\sim (-\log x)^{-\gamma}\ell(-\log x) \ \ \text{as}
\ x\to 0,
\end{equation}
holds with
 $\ell$ slowly varying at $\infty$ and
$\gamma\in [1,2)$, and assume that ${\tt m}<\infty$ if $\gamma=1$, then for
$b_n={\tt m}^{-1}\log n$, $a_n={\tt m}^{-(\gamma+1)/\gamma}c_{\lfloor\log
n\rfloor}$  and $c_n$  any sequence satisfying
$\lim_{n\to\infty}n\ell(c_n)/c_n^\gamma=1$, the limiting distribution
is $\gamma$-stable with characteristic function
$$\tau\mapsto \exp\{-|\tau|^\gamma
\Gamma(1-\gamma)(\cos(\pi\gamma/2)+i\sin(\pi\gamma/2)\, {\rm
sgn}(\tau))\}, \ \tau\in{\mathbb R}.$$

\item[{\rm (d)}] Assume that ${\tt m}=\infty$ and the
relation {\rm (\ref{domain1})} holds with $\gamma=1$. Let $c$ be any
positive function satisfying $\lim_{x\to\infty}x\ell(c(x))/c(x)=1$
and set 
$$\psi(x):=x\int_{\exp(-c(x))}^1 \prob(1-W\leq y)/y\,{\rm d}y.$$
Let $b$ be any positive function satisfying
$b(\psi(x))\sim\psi(b(x))\sim x$ 
(asymptotic inverse to $\psi$). 
Then, with $b_n=b(\log n)$ and
$a_n=b(\log n)c(b(\log n))/\log n$, the limiting distribution is
$1$-stable with characteristic function
\begin{equation}\label{stable}
\tau\mapsto \exp\{-|\tau|(\pi/2-i\log|\tau|\,{\rm sgn}(\tau))\}, \ \tau\in{\mathbb R}.
\end{equation}

\item[{\rm (e)}] If the relation {\rm (\ref{domain1})} holds with $\gamma \in [0,1)$
then, for $b_n=0$ and $a_n:=\log^\gamma n/\ell(\log n)$, the
limiting distribution  is the scaled Mittag-Leffler law $\theta_\gamma$
(exponential, if $\gamma=0$) characterised by the moments
$$\int_0^\infty
x^n\theta_\gamma({\rm d}x)=\dfrac{n!}{\Gamma^n(1-\gamma)\Gamma(1+n\gamma)}\,\,,
~~~ n\in\Nat.$$
\end{enumerate}
\end{theorem}
\noindent
\proof
The results are first derived for $M_n$ by adopting asymptotics from 
the renewal theory. 
To pass to $K_n$ it is shown, 
under the condition ${\tt m}_1<\infty$, that the variable  $M_n-K_n$ (the number 
of empty boxes to the left of $U_{n:n}$) converges in distribution and in the mean to a random variable 
with expected value ${\tt m}_1/{\tt m}$.

\endpf

\vskip0.2cm
\noindent
{\bf Example} Suppose $W$ has a beta density (\ref{beta}).
The moments are easily computable as 
${\tt m}=\Psi(\theta+\gamma)-\Psi(\theta),~~{\tt m}_1=\Psi(\theta+\gamma)-\Psi(\gamma),~~\sigma^2=\Psi'(\theta)-\Psi'(\theta+\gamma)$
(with $\Psi=\Gamma'/\Gamma$ denoting the logarithmic derivative of the gamma function).
We are therefore in the case (a) of Theorem \ref{main1}, hence $K_n$ is asymptotically normal
with $\ex[K_n]\sim {\tt m}^{-1}\log n$ and ${\rm Var}[K_n]$ of the same order.

\par The instance $\gamma=1$ recovers well-known asymptotics of Ewens' partitions, 
which have $K_n\sim\theta\log n$ a.s.
In this case the limit law of
the number of empty gaps
$M_n-K_n$ 
has probability generating function $z\mapsto \Gamma(\theta+1)\Gamma(\theta+1-z\theta)/\Gamma(1+2\theta-\theta z)$
(which identifies a mixture of Poisson distributions, see \cite{GIksetal}).

\vskip0.2cm
\noindent
{\bf Example} Suppose the law of $W$ is given by $\prob(1-W\leq x)=(1-\log x)^{-1}$, $x\in (0,1)$.
It can be checked that ${\tt m}_1<\infty$, hence the case (c) applies and 
$${(\log \log n)^2\over \log n}K_n-\log \log n-\log\log\log n$$
converges to a $1$-stable law with characteristic function (\ref{stable}).
The number of empty boxes $M_n-K_n$ converges in probability to $0$.

\vskip0.2cm

\par Under assumptions ${\tt m}<\infty, {\tt m}_1<\infty$ the limit behaviour of $K_{n,r}$'s
is read from a limiting occupancy model \cite{Small}. 
 To describe the limit
we 
pass 
to the dual composition, generated by right-to-left stick-breaking
${\cal R}=\{V_1\cdots V_i: i\geq 1\}$ with independent $1-V_i\ed W$. Let $(X_{n,1},X_{n,2},\ldots)$ be the occupancy numbers
of the gaps in the left-to-right order, this is a random weak composition ($0$'s allowed) of $n$ 
with $X_{n,1}>0$ and $X_{n,j}\geq 0$ for $j>1$.
By inflating $[0,1]$ with factor $n$, the uniform sample converges as a point process to a unit Poisson process (balls).
On the other hand, 
$n{\cal R}$ converges to a self-similar point process ${\cal Z}$, whose gaps play the role of boxes.
 From this, the occupancy vector $(X_{n,1},X_{n,2},\ldots)$ acquires a limit, which is  
an occupancy vector $(X_1,X_2,\ldots)$ derived from the limiting point processes. The limit distribution  of the occupancy vector
is
$$\prob(X_1=\lambda_1,\ldots,X_k=\lambda_k)= {1\over {\tt m}(\lambda_1+\ldots+\lambda_k)}\prod_{i=1}^k \hat{q}(\Lambda_i:\lambda_i)$$
where $\lambda_1>0,\lambda_i\geq 0,\Lambda_i=\lambda_1+\ldots+\lambda_i$ and  $\hat{q}(n:m)={n\choose m}\ex[W^m(1-W)^{n-m}]$.
Correspondingly,
$K_{n,r}$'s jointly converge in distribution to $\#\{i:X_i=r\}$, $r=1,2,\ldots$.
The convergence also holds for $K_{n,0}$, defined as the number of g4

\par If $W\ed{\rm beta}(1,\theta)$ then $\cal Z$ is Poisson process with density $\theta/x$.
Then
$K_{n,r}$'s  converge in distribution  to independent Poisson variables with mean $\theta/r$, 
which is a well known property of Ewens' partitions 
\cite{ABT}.
It is a challenging open problem to identify 
the limit laws of the $K_{n,r}$'s for general distribution of $W$.

\subsection{Regular variation: $0<\alpha\leq 1$.}

Suppose (\ref{regvar}) holds with $0<\alpha\leq 1$.
This case 
is treated by reducing the occupancy problem to counting the gaps of given sizes.
For $x>0$ let $N_x(t)$ be the number of gaps of size at least $x$, in the partial range of the multiplicative subordinator
 $\big(1-\exp(-S_u), 0\leq u<t\big)$.     
Introduce the exponential functionals
$$
I_\alpha(t):=\int_0^t \exp(-\alpha S_t)\,\dd t,~~~~I_\alpha:=I_\alpha(\infty).
$$
The distribution of $I_\alpha(\infty)$ is determined by the formula for the moments  \cite{Carmona}
$$
\ex (I_\alpha)^k={k!\over\prod_{i=1}^k\Phi(\alpha j)}\,,
$$
where $\Phi$ is the Laplace exponent of the subordinator $(S_t)$.

\begin{theorem}\label{gaps}{\rm \cite{GPYII}}
Suppose the L{\'e}vy measure fulfills {\rm (\ref{regvar})}.
Then for $0<t\leq \infty$
\begin{eqnarray*}
 {\rm for~~} 0<\alpha<1~~~~~{N_x(t)\over \ell(1/x)x^{-\alpha}}&\to& I_\alpha(t)~~{\rm a.s.},~~x\darr 0,\\
{\rm for~~} \alpha=1~~~~~{N_x(t)\over \ell_1(1/x)x^{-\alpha}}& \to & I_\alpha(t)~~{\rm a.s.},~~x\darr 0,
\end{eqnarray*}
where 
$$\ell_1(z)=\int_z^\infty u^{-1}\ell(u) \dd u$$
is another function of slow variation, satisfying $\ell_1(z)\gg \ell(z)$ as $z\to\infty$.
\end{theorem}
\proof
Let $\widetilde{N}_{x}(t)$ be the number of gaps in the range of the (additive) subordinator restricted to $[0,t]$.
By the L{\'e}vy-Ito construction of $(S_t)$ from a Poisson process, we have the strong law 
$\widetilde{N}_{x}(t)\sim \nut[y,\infty] t$ a.s. for $y\darr 0$. 
A small gap $(s,s+x)$  is mapped by the function $s\to 1-e^{-s}$ in a gap of size $e^{-s}x$, from which the result 
for finite $t$ follows by integration.
Special tail estimates are required to conclude that similar  asymptotics hold with integration  extended to $[0,\infty]$.
\endpf
\noindent
The instance $\alpha=1$ may be called in this context the case of {\it rapid variation}. In this case
$\ell$ in (\ref{regvar}) must decay at $\infty$ sufficiently fast, in order to satisfy $\Phi(1)<\infty$.

\par 
Conditioning on the frequencies ${\cal S}=(s_j)$ embeds the small-block problem in the framework of the classical 
occupancy problem:
$n$ balls are thrown 
in an infinite series of boxes, with positive probability $s_j$ of hitting box $j$.
By a result of Karlin \cite{Karlin}, the number of occupied boxes is  asymptotic to the expected 
number, from which 
$K_n\sim  \ex[K_n\giv{\cal R}]$ a.s., and a similar result holds for 
$K_{n,r}$ under the regular variation with index $\alpha>0$. Combining this with Theorem \ref{gaps}, we have
(see  \cite{GPYII})

\begin{theorem}
\label{asympRV}
Suppose the L{\'e}vy measure fulfills {\rm (\ref{regvar})}. Then,  uniformly in $0<t\leq\infty$, as $n\to\infty$, 
the convergence holds almost surely and in the mean:
$$
{\KK_n(t)\over \Gamma(1-\alpha)n^\alpha\ell(n)}\to I_\alpha(t), ~~~~
{\KK_{n,r}(t)\over \Gamma(1-\alpha)n^\alpha\ell(n)}\to (-1)^{r-1}{\alpha\choose r}I_\alpha(t), 
$$
for $0<\alpha<1$ and $r\geq 1$, or $\alpha=1$ and $r>1$. Similarly, $\KK_n(t)/(n\ell_1(n))\to I_1(t)$ for $\alpha=1$.
\end{theorem} 
\noindent
Thus $K_n, K_{n,r}$ have the same order of growth if $0<\alpha<1$. In the case $\alpha=1$ of rapid variation,
singletons dominate, $K_{n,1}\sim K_n$, while all other $K_{n,r}$'s with $r>1$ are of the same order of growth
which is smaller than that of $K_n$. 

\vskip0.2cm
\noindent
{\bf Example} The subordinator associated with the two-parameter family of compositions has $\Phi$ given by (\ref{althlap}), hence
$$
\ex (I_\alpha)^k={(\alpha+\theta)(2\alpha+\theta)\cdots((k-1)\alpha+\theta)\Gamma(\theta+1)\over \Gamma(k\alpha+\theta)\{\alpha\Gamma(1-\alpha)\}^k},
$$
which for $\theta=0$ and $0<\alpha<1$ identifies the law of $I_\alpha$ as a Mittag-Leffler distribution.
Theorem \ref{asympRV} recovers in this instance known asymptotics \cite{CSP}  related 
to the local times of Bessel bridges.

\par One generalisation of Theorems \ref{gaps} and \ref{asympRV} is obtained by taking for paintbox the range of a process
$\phi(S_t)$, where $\phi:{\mathbb R}_+\to [0,1]$ is a smooth monotone function, with not too bad behaviour at $\infty$.
The generalised power laws hold also for this class of partitions, 
with the only difference that
the exponential functionals should be replaced  by  integrals $\int_0^t \{\phi'(S_u)\}^\alpha\dd u$, see \cite{GPYII}.

\subsection{Slow variation: $\alpha=0$}

The case of infinite L{\'e}vy measure with slowly varying tail $\nut[y,\infty]\sim \ell(1/y)$ $(y\darr 0)$ is intermediate between
finite $\nut$ and the case of proper regular variation.  
In this case  $K_{n,r}\to\infty$ (like in the case $\alpha>0$) but $K_{n,r}\ll K_n$ (like in the case of finite $\nut$). 
Following 
Barbour and Gnedin \cite{BG} 
we will exhibit
a further 
wealth of possible modes of asymptotic behaviour
appearing in this transitional regime. 

\par 
We assume  that the first two moments of the subordinator are finite.
The assumption about the moments is analogous to the instance (a) of Theorem \ref{main1} in the case of finite $\nut$. 
The results will be formulated for the case
$$ \ex[ S_t]=t,~~~~~\Var[ S_t]={\tt s}^2 t,$$
which can be always achieved by a linear time scaling.
Indeed, a general subordinator $S_t$ with 
$${\tt m}:=\ex [S_1]=\int_0^\infty x\nut(\dd x),~~~v^2:=\Var[S_1]=\int_0^\infty x^2\nut(\dd x)$$
should be replaced by $S_{t/{\tt m}}$, then ${\tt s}^2=v^2/{\tt m}$. 
Because the linear time change 
does not affect the range of the process, it does not change the distribution of 
$K_n, K_{n,r}$.

\par For the sample (balls) we take a Poisson point process on $[0,1]$ with intensity  $n>0$.
This is the same as assuming a Poisson$(n)$ number of uniform points thrown on $[0,1]$.
To avoid new notations, we
further understand $n$ as the intensity parameter, 
and use the old notation 
$\KK_n(t)$ to denote  the  number of blocks of the 
(poissonised) subcomposition on the interval  $[0, 1-\exp(- S_t)]$. The convention  for $K_n$ is the same. 
For large samples the poissonised quantities are very close to   
 their fixed-$n$ counterparts, but the Poisson framework is easier to work with.

\par The total number of blocks is the terminal value $K_n=\KK_n(\infty)$ of the increasing process $\KK_n(t)$.
Poissonisation makes the subcompositions within $[0,1-\exp(-S_t)]$ and $[1-\exp(-S_t),\,1]$ conditionally independent given $S_t$,
hence $\KK_n(t)$ and $\KK_n(\infty)-\KK_n(t)$ are also conditionally independent.
The consideration can be restricted to the time range $t<\tau_n$ where 
$\tau_n:=\inf\{t:S_t>\log n\}$ is the passage time through $\log n$, since after this time the number of blocks produced is bounded
by a Poisson(1) variable.

\par Define the poissonised Laplace exponent
$$
\Phi_0(n):=\int_0^\infty \{1-\exp(-n(1-e^{-y}))\}\nut(\dd y).
$$
For large $n$, we have $\Phi_0(n)\sim \Phi(n)$, but the former is more convenient to deal with, 
since it enters naturally the {\it compensator} of $(\KK_n(t), t\geq 0)$,
$$
A_n(t):=\int_0^t \Phi_0(n\exp(-S_u))\dd u.
$$

\par Introduce 
$$\Phi_k(n):=\int_0^\infty \{\Phi_0(ne^{-s})\}^k \dd s= \int_0^n \{\Phi_0(s)\}^k\,{\dd s\over s}\,,~~~~~~k=1,2.$$
By the assumption of slow variation and from $\Phi_0(n)\to \infty$  it follows that $\Phi_k$'s are also slowly varying, and satisfy
$$\Phi_2(n)\gg\Phi_1(n)\gg\Phi_0(n),~~~~n\to\infty.$$
These functions give, asymptotically, the moments of $K_n$ and of the terminal value of the compensator
$$\ex[K_n]=\ex[A_n(\infty)]\sim \Phi_1(n),~~\Var[K_n]\sim\Var[A_n(\infty)]\sim{\tt s}^2\Phi_2(n),~~~~~n\to\infty.$$
\vskip0.2cm
\noindent
{\bf Remark} In the stick-breaking case $\nut[0,\infty]=1$ the asymptotics of $\Var[A_n(\infty)]$ and $\Var[K_n]$ are different,
because the asymptotic relation  $\Phi_1(n)\ll\Phi_2(n)$ is not valid. 
Instead, we have $\Var[A_n(\infty)]\sim v^2{\tt m}^{-3}\log n$, and ~$\Var[K_n]\sim 
\sigma^2{\tt m}^3\log n$ with $\sigma^2=v^2-{\tt m}^2$.

\vskip0.2cm

\par The following approximation lemma reduces the study of $\KK_n(t)$ to the asymptotics of the compensator.

\begin{lemma}
We have, as $n\to\infty$,
$$
\ex[K_n-A_n(\infty)]^2\sim\Phi_1(n),~~~
$$
and for any $b_n$ such that $\Phi_1(n)/b_n^2\to 0$
$$
\lim_{n\to\infty} \prob\left[\sup_{0\leq t\leq\infty}|\KK_n(t)-A_n(t)|>b_n
\right]=0.
$$
\end{lemma}  
\proof Noting that 
$\KK_n(t)-A_n(t)$ is a square integrable martingale with unit jumps, we derive $\ex[K_n-A_n(\infty)]^2=\ex[A_n(\infty)]$,
from which the first claim follows. The second follows by application of Kolmogorov's inequality.
\endpf

\noindent
From this the law of large number is derived:
\begin{theorem} As $n\to\infty$, we have 
$K_n\sim A_n(\infty)\sim\Phi_1(n)$ 
almost surely and in the mean.
\end{theorem}

\par For more delicate results we need to
 keep fluctuations of the compensator under control
For this purpose we adopt one further assumption,
akin to  de Haan's second order regular variation \cite{BGT}.
As in  Karamata's representation of slowly varying functions  \cite{BGT},
write $\Phi_0$ as
$$
\Phi_0(s)=\Phi_0(1)\exp\left(\int_0^s {\dd z\over z L(z)} \right),
$$
where 
$$
L(n):={\Phi_0(n)\over n\Phi_0'(n)}\,.
$$

\vskip0.2cm
\noindent
{\bf The key assumption.} There exist constants $c_0, n_0>0$ such that
\eq\label{key-a}
\left|
{nL'(n)\over L(n)}\right| <{c_0\over\log n}\,\,~~~{\rm for~}n>n_0.
\en
\vskip0.2cm
\noindent
In particular, $L$ is itself slowly varying, which is equivalent to the slow variation of $n\Phi_0'(n)$ as $n\to\infty$.
Note that the faster $L$, the slower $\Phi_0$.
The assumption allows to limit local variations of $\Phi_0$, 
which makes possible
approximating
the compensator by a simpler process 
$$
A_n^*(t):=\int_0^{t\wedge\log n} \Phi_0(ne^{-u})\left(1- {S_u-u\over L(ne^{-u})}\right) \dd u,
$$
in which the subordinator enters linearly.
This in turn allows to derive the limit behaviour of the compensator from the functional CLT for $(S_t)$ itself.

\subsubsection {Moderate growth case} \label{moderate}

This is the case  $L(n)\asymp\log n$.
We shall state weak convergence in the space
$D_0({\mathbb R}_+)$ of c{\'a}dl{\'a}g 
functions with finite limits at $\infty$.

\par The time-changed scaled  process  
$$
\KK^{(1)}_n(u):=
(\Phi_0(n)\sqrt{\log n})^{-1}\left( \KK_n(u\log n)-\log n\int_0^{u\wedge 1}\Phi_0(n^{1-v})\dd v\right)
$$
converges weakly to the process
$$
Y_n^{(1)}(u):={\tt s}\int_0^{u\wedge 1} h_n^{(1)}(v)B_v\dd v,
$$
where $(B_u)$ is the BM and
$$
h_n^{(1)}(u):={\Phi_0(n^{1-u})\log n\over \Phi_0(n)L(n^{1-u})}\,. 
$$
 In particular, if $L(n)\sim\gamma \log n$ for some $\gamma>0$, we have 
$$
Y^{(1)}(u)={\tt s} \int_0^{u\wedge 1} \gamma^{-1}(1-v)^{(1-\gamma)/\gamma} B_v \dd v.
$$
\vskip0.2cm
\noindent
{\bf Example} Consider a subordinator with Laplace exponent $\Phi(n)\sim c\log^{1/\gamma} n$, 
${\tt m}=\ex[S_1]=\Phi'(0)  ,~{v}^2=\Var[S_1]=\Phi''(0)$.
A CLT for $K_n$ holds with standard scaling and centering by the moments
$$
\ex[K_n]\sim  {c\log^{1+1/\gamma} n\over{\tt m}(1+1/\gamma)}\,,~~~~\Var[K_n]\sim {c^2{v}^2\log^{1+2/\gamma} n\over {\tt m}^3 (1+2/\gamma)}.
$$ 
 
A special case is the gamma subordinator with 
$$\nut(\dd y)= \theta e^{-\theta y} \dd y/y, ~\Phi(n)=\theta\log(1+n/\theta),~{\tt m}=1,
~{v}^2=1/\theta.$$
Some generalisations are considered in \cite{GPYI}.

\subsubsection{Fast growth case}

This case is defined by the conditions  $L(n)\to\infty,~L(n)\ll \log n$,
then $\Phi_0$ grows faster than any power of the logarithm. For instance 
$\Phi_0(n)\asymp \exp(\log^\gamma n)$ with $0<\gamma<1$.
The scaled process 
$$
\KK_n^{(2)}(u):=(\Phi_0(n)\sqrt{L(n)})^{-1}\left(
\KK_n(uL(n))-L(n)\int_0^{u\wedge(\log n/L(n))} \Phi_0(n\exp(-v L(n))\dd v
\right)
$$
converges weakly to
$$
Y^{(2)}(u):= {\tt s}\int_0^u e^{-v}B_v \dd v.
$$

\subsubsection{Slow growth case}

Suppose that $L(n)=c(n) \log n$, where $c(n)\to\infty$ but slowly enough to have  
$$\int_2^\infty {\dd n\over c(n) n \log n}=\infty $$
(otherwise $\nut$ is a finite measure). 
For instance we can have $c(n)\asymp \log\log n$ (in which case $\Phi(n)\asymp\log\log n$),
but the growth $c(n)\asymp\log^\gamma n$ with $\gamma>0$ is excluded.
Like in the case of finite $\nut$, almost all variability of $K_n$ 
comes from the range of times very close to the passage time $\tau_n$.

\par 
The key quantity describing the process $ \KK_n(t)$ in this case is the family of integrals
$$\int_{(\tau_n-t)_+}^{\tau_n} \Phi_0(e^v)\dd v, ~~~t\geq 0,$$
where $\tau_n$ is the passage time at level $\log n$.
The randomness enters here only through $\tau_n$, which is approximately normal for large $n$ with $\ex[\tau_n]\sim \log n$,
$\Var[\tau_n]\sim {\rm s}^2\log n$.
The process
$$
\KK_n^{(3)}(t):=(\Phi_0(n)\sqrt{\log n})^{-1} \left(\KK_n(t)-\int_{(\log n-t)_+}^{\log n} \Phi_0(e^v)\dd v
\right)
$$
is approximated by 
$$
Y^{(3)}(u)=:{\tt s}\eta-(\Phi_0(n)\sqrt{\log n})^{-1} \int_{(\log n-t)_+}^{(\log n -t+{\tt s}\eta\sqrt{\log n})_+} \Phi_0(e^v)\dd v,
$$
where $\eta$ is a standard normal random variable.

\subsubsection{Gamma subordinators and the like}

Asymptotics of $K_{n,r}$ is known \cite{GPYI} for gamma-like subordinators with logarithmic behaviour 
$\nut[y,\infty]\sim -c\log y$ for $y\to 0$, under some additional assumptions on the tail of $\nut$ for $y$ near $0$ and $\infty$.
This case  is well suited for application of singularities analysis to formulas like
$$
\int_0^\infty n^{s-1}\phi_r(n)\dd n= {\Gamma(r+s)\over r!} \,{\Phi(-s:-s)\over\Phi(-s)}\,,~~~~-1<\Re s<0
$$
for the Mellin transform of the expected value $\phi_r(n)$ of poissonised $K_{n,r}$. In this formula $\Phi(s)$ and $\Phi(-s:-s)$ 
are the analytical continuations  in the complex domain
of the Laplace exponent $\Phi(n)$ and the bivariate funtion $\Phi(n:n)$, respectively.

For the moments we have 
\begin{eqnarray*}
\ex[K_n]&\sim & {\log^2 n\over 2 {\tt m}},~~\Var[K_n] \sim  {v^2\log^3 n\over 3{\tt m}},\\
\ex[K_{n,r}]&\sim&{\log n\over r{\tt m}},~~\Var[K_{n,r}]\sim  \left({v^2\over r^2 {\tt m}^3}+{1\over r{\tt m}}\right)\log n.
\end{eqnarray*}
The CLT for $K_n$ is an instance of the moderate growth case in Section \ref{moderate}.

\par As $n\to\infty$, the infinite sequence of scaled and centered block counts
$$\left({K_{n,r}-\ex[K_{n,r}]\over \sqrt{\log n}}\,,r=1,2,\ldots\right)$$ 
converges in distribution to a multivariate Gaussian sequence with the covariance matrix
$$
{v^2\over {\tt m}^3}\,{1\over ij}+1(i=j) {1\over j{\tt m}}\,,~~~~~~i,j=1,2,\ldots.
$$
See \cite{GPYI} for explicit assumptions on $\nut$ in this logarithmic case  and further examples, including computations for the subordinator in 
Example \ref{gammageom}.
\vskip0.2cm
\par The behaviour of $K_{n,r}$'s for other slowly varying infinite $\nut$ remains an open problem.

\vskip0.5cm
\noindent
{\bf Acknowledgement} 
These notes were presented in a series of lectures at the School on Information and Randomness 2008, in Santiago de Chile.
The author is indebted to Servet Martinez for the invitation and motivation to record the lectures.

\end{document}